\documentclass[11pt]{article}
\usepackage{amssymb,amsmath,latexsym}

\begin{document}

\begin{doublespace}

\newtheorem{thm}{Theorem}[section]
\newtheorem{lemma}[thm]{Lemma}
\newtheorem{defn}[thm]{Definition}
\newtheorem{prop}[thm]{Proposition}
\newtheorem{corollary}[thm]{Corollary}
\newtheorem{remark}[thm]{Remark}
\newtheorem{example}[thm]{Example}
\numberwithin{equation}{section}
\def\ee{\varepsilon}
\def\qed{{\hfill $\Box$ \bigskip}}
\def\NN{{\cal N}}
\def\AA{{\cal A}}
\def\MM{{\cal M}}
\def\BB{{\cal B}}
\def\LL{{\cal L}}
\def\DD{{\cal D}}
\def\FF{{\cal F}}
\def\EE{{\cal E}}
\def\QQ{{\cal Q}}
\def\RR{{\mathbb R}}
\def\R{{\mathbb R}}
\def\L{{\bf L}}
\def\E{{\mathbb E}}
\def\F{{\bf F}}
\def\P{{\mathbb P}}
\def\N{{\mathbb N}}
\def\eps{\varepsilon}
\def\wh{\widehat}
\def\pf{\noindent{\bf Proof.} }

\title{\Large \bf  Potential Theory of Truncated Stable Processes}
\author{Panki Kim\thanks{The research of this author is supported by Research Settlement Fund for the new faculty of SNU.}\\
Department of Mathematics\\
Seoul National University\\
Seoul 151-742, Republic of Korea\\
Email: pkim@snu.ac.kr \smallskip \\
and
\smallskip \\
Renming Song\thanks{The research of this author is supported in part
by a joint
US-Croatia grant INT 0302167.}\\
Department of Mathematics\\
University of Illinois \\
Urbana, IL 61801, USA\\
Email: rsong@math.uiuc.edu }
\date{ }
\maketitle

\begin{abstract}
For any $\alpha\in (0, 2)$, a truncated symmetric $\alpha$-stable
process is a symmetric L\'evy process in $\R^d$ with a  L\'evy
density given by $c|x|^{-d-\alpha}\, 1_{\{|x|< 1\}}$ for some
constant $c$. In this paper we study the potential theory of
truncated symmetric stable processes in detail. We prove a Harnack
inequality for nonnegative harmonic nonnegative functions of these
processes. We also establish a boundary Harnack principle for
nonnegative functions which are harmonic with respect to these
processes in bounded convex domains. We give an example of a
non-convex domain for which the boundary Harnack principle fails.
\end{abstract}



\noindent {\bf Running Title:} Truncated Stable Processes

\section{Introduction}

Recently there have been a lot of interests in studying
discontinuous stable processes due to their importance in theory as
well as applications. Many important results have been established.
These results include, among other things, sharp estimates on the
Green functions and Poisson kernels (\cite{CS1} and \cite{Ku1}), the
boundary Harnack principle (\cite{B} and \cite{SW}) and the
identification of the Martin boundary for various domains
(\cite{B2}, \cite{CS} and \cite{SW}). See \cite{C1} for a survey of
some of these results.

However in a lot of applications one needs to use discontinuous
Markov processes which are not stable processes. Therefore we need
to extend the known results on stable processes to other discontinuous
Markov processes.

In \cite{R} and \cite{CS3}, sharp estimates on the Green functions
of killed relativistic stable processes in bounded $C^{1, 1}$ domains
were established. These estimates can be used to establish various
properties of relativistic stable processes.

Another discontinuous Markov process,
the censored stable process,
was introduced and studied in \cite{BBC}.
Roughly speaking, for $\alpha \in (0, \, 2)$,
a censored $\alpha$-stable process
in an open set  $D\subset \R^d$
is a process obtained from a symmetric $\alpha$-stable
L\'evy process
by restricting its L\'evy measure
to $D$.
The censored process is repelled from the
complement of the open set $D$ because it is prohibited to
make jumps outside $D$. Some potential theoretic properties
of censored stable processes, such as Green function
estimates, Martin boundary, and Fatou type theorem,
were established  recently
(see \cite{CK1}, \cite{CK2} and \cite{K}).

In this paper we study yet another type of discontinuous Markov
processes which we call truncated symmetric stable processes. For
$\alpha\in (0, 2)$, a truncated symmetric $\alpha$-stable process is
a symmetric L\'evy process in $\R^d$ whose L\'evy density $l(x)$
coincides with the L\'evy density of a symmetric $\alpha$-stable
process for $|x|$ small (say, $|x|<1$) and is equal to zero for
$|x|$ large (say, $|x|\ge 1$). In other words, a truncated symmetric
$\alpha$-stable process is a symmetric L\'evy process in $\R^d$ with
a  L\'evy density given by $c|x|^{-d-\alpha}\, 1_{\{|x|< 1\}}$ for
some constant $c$. Truncated stable processes are very natural and
important in applications where only jumps up to a certain size are
allowed. One expects that many properties of the truncated stable
processes should be similar to those of the symmetric stable
processes, but some properties are very different. For instance, the
boundary Harnack principle for symmetric stable processes is valid
on any $\kappa$-fat set, while we will show that on non-convex
domains the boundary Harnack principle for truncated stable
processes might fail.

In some aspects, truncated stable processes have
nicer behaviors and are more preferable than symmetric
stable processes, for instance, by Theorem 25.17 of \cite{Sa}
we know that truncated stable processes have finite exponential
moments. However, as we shall see later, in some other respects,
truncated stable processes
are much more difficult and more delicate to study than symmetric
stable processes.

The starting point of our research on truncated stable processes was
our attempt to establish a Harnack inequality for nonnegative
harmonic functions of truncated stable processes. The recent
developments in Harnack inequalities for discontinuous Markov
processes were initiated in \cite{BL}. The method of \cite{BL} was
extended in \cite{SV} to cover a large class of Markov processes.
Two other methods for proving the Harnack inequality for
discontinuous Markov processes were contained in \cite{BSS} and
\cite{CK}. However, none of the methods above apply to truncated
stable processes. This gives another indication that truncated
stable processes are pretty delicate to deal with.

Our strategy for studying truncated stable processes is as follows.
First, we consider
killed truncated stable processes on small sets and show that its
Green functions are comparable to
the Green functions of the corresponding
killed symmetric stable processes. Then we study Poisson kernels
for truncated stable processes on small sets in detail. Finally
we prove the Harnack inequality
and boundary Harnack principle for nonnegative harmonic
functions of truncated stable processes by using properties of its
Poisson kernels and some ideas in
\cite{B}, \cite{BB} and \cite{SW}.

In this paper we will always assume that $d\ge 2$. The case of $d=1$
can also be considered, but some arguments need to be modified. We
leave this case to the interested reader.

In this paper, we use ``$:=$" as a way of
definition, which is  read as ``is defined to be".
The letter $c$, with or
without subscripts, signifies a constant whose value is
unimportant and which may change from location to location,
even within a line.

\section{Stable Processes and Truncated Stable Processes}

Throughout this paper we assume  $\alpha\in (0, 2)$ and $d \ge 2$.
Recall that a symmetric $\alpha$-stable process $X=(X_t, \P_x)$ in $\R^d$
is a L\'evy process such that
$$
\E_x\left[e^{i\xi\cdot(X_t-X_0)}\right]=e^{-t|\xi|^{\alpha}},
\quad \quad \mbox{ for every } x\in \R^d \mbox{ and } \xi\in \R^d.
$$
The Dirichlet form $({\cal E}, D({\cal E}))$
associated with $X$ is given by
$$
{\cal E}(u, v):=\int_{\R^d} \hat{u} (\xi)\bar{\hat{v}}(\xi)|\xi|^{\alpha}
d\xi,
\quad
D({\cal E}):=
\{ u\in L^2(\R^d):\int_{\R^d} |\hat{u} (\xi)|^2|\xi|^{\alpha}d\xi
<\infty\},
$$
where
$\hat{u} (\xi):=(2\pi)^{-d/2}\int_{\R^d}e^{i\xi\cdot y}u(y)dy$ is the
Fourier transform of $u$.
As usual, we define
${\cal E}_1(u, v):={\cal E}(u, v)+(u, v)_{L^2(\R^d)}$ for
$u, v\in D({\cal E})$.
Then we have
\begin{equation}\label{1}
{\cal E}_1(u, u)=\int_{\R^d} |\hat{u} (\xi)|^2(1+|\xi|^{\alpha})d\xi,
\quad u\in D({\cal E}).
\end{equation}
Another expression for ${\cal E}$ is as follows:
$$
{\cal E}(u, v)=\frac12
{\cal A} (d, - \alpha)\int_{\R^d}\int_{\R^d}
\frac{(u(x)-u(y))(v(x)-v(y))}
{|x-y|^{d+\alpha}}\, dxdy,
$$
where
$
{\cal A}(d, -\alpha):= \alpha2^{\alpha-1}\pi^{-d/2}
\Gamma(\frac{d+\alpha}2)
\Gamma(1-\frac{\alpha}2)^{-1}.
$
Here $\Gamma$ is the Gamma function defined by
$\Gamma(\lambda):= \int^{\infty}_0 t^{\lambda-1} e^{-t}dt$
 for every $\lambda > 0$.
By a truncated symmetric $\alpha$-stable process in $\R^d$ we mean a
symmetric L\'evy
process $Y=(Y_t, \P_x)$ in $\R^d$
such that
$$
\E_x\left[e^{i\xi\cdot(Y_t-Y_0)}\right]=e^{-t\psi(\xi)},
\quad \quad \mbox{ for every } x\in \R^d \mbox{ and } \xi\in \R^d,
$$
with
\begin{equation}\label{e:psi}
\psi(\xi)={\cal A} (d, - \alpha) \int_{\{|y|<
1\}}\frac{1-\cos(\xi\cdot y)}{|y|^{d+\alpha}}dy.
\end{equation}
The Dirichlet form $({\cal Q},
D({\cal Q}))$ of $Y$ is given by
$$
{\cal Q}(u, v):=\int_{\R^d} \hat{v} (\xi)\bar{\hat{u}}(\xi)\psi(\xi)
d\xi,
\quad
D({\cal Q}):=
\{ u\in L^2(\R^d):\int_{\R^d} |\hat{u} (\xi)|^2\psi(\xi)d\xi
<\infty\}.
$$
The Dirichlet form
$({\cal Q}, D({\cal Q}))$ of $Y$ can also be written as follows
\begin{eqnarray*}
{\cal Q}(u, v)&=&
\frac12
{\cal A} (d, - \alpha)\int_{\R^d}\int_{\R^d}
\frac{(u(x)-u(y))(v(x)-v(y))}
{|x-y|^{d+\alpha}}1_{\{|x-y|< 1\}}dxdy\\
D({\cal Q})&=&
\left\{ u\in L^2(\R^d):\int_{\R^d}\int_{\R^d}
\frac{(u(x)-u(y))^2}
{|x-y|^{d+\alpha}}1_{\{|x-y|< 1\}}dxdy<\infty\right\}.
\end{eqnarray*}
Similar to ${\cal E}_1$, we can also define ${\cal Q}_1$.
Then we have
\begin{equation}\label{2}
{\cal Q}_1(u, u)=\int_{\R^d} |\hat{u} (\xi)|^2(1+\psi(\xi))
d\xi.
\end{equation}
By the change of variable $y=x/|\xi|$, we have from (\ref{e:psi})
$$
\psi(\xi)={\cal A} (d, - \alpha)|\xi|^{\alpha}
\int_{\{|x|<|\xi|\}}\frac{1-\cos(\frac{\xi}{|\xi|}\cdot x)}
{|x|^{d+\alpha}}dx.
$$
Since $1-\cos(\frac{\xi}{|\xi|}\cdot x)$ behaves like $|x|^2$
for small $|x|$,
it is easy to check that $\psi(\xi)$ behaves like $|\xi|^2$ near the origin.
Also we see that as  $|\xi|$ goes to infinity, the integral in
the above equation
goes to a positive constant. So  $\psi(\xi)$
 behaves like $|\xi|^{\alpha}$ near infinity. Therefore
 by the definition
of $D({\cal E})$ and $D({\cal Q})$ we see that $D({\cal E})=D({\cal Q})$.
From now on we will use ${\cal F}$ to stand for $D({\cal E})$. Using
(\ref{1}), (\ref{2}) and the fact above, we see that there exist positive
constants $c_1$ and $c_2$ such that
$$
c_1 {\cal E}_1(u, u)\le {\cal Q}_1(u, u)\le c_2 {\cal E}_1(u, u),
\quad u\in {\cal F}.
$$
Therefore the capacities corresponding to the Dirichlet forms
$({\cal E}, {\cal F})$ and $({\cal Q}, {\cal F})$ are comparable,
hence we get that a set $A$ has zero capacity with respect to
$({\cal E}, {\cal F})$ if and only if it has
zero capacity with respect to
$({\cal Q}, {\cal F})$ and that
a function $u$ is quasi continuous with respect to the
capacity of $({\cal E}, {\cal F})$ if and only if it is
quasi continuous with respect to the
capacity of $({\cal Q}, {\cal F})$.
So when we speak of quasi continuous functions or sets of zero capacity,
we do not need to specify the Dirichlet forms.
For concepts and results related to Dirichlet forms, we refer our readers to
\cite{FOT}.

It is well known that any function $u\in {\cal F}$ admits a
quasi continuous version. From now on, whenever we talk about
a function $u\in {\cal F}$, we always use the quasi continuous version.

Using the asymptotic behavior of $\psi$ and Proposition 28.1 in
\cite{Sa} we know that the process $Y$ has a smooth density $p^Y(t,
x, y)$. Since $\psi(\xi)$ behaves like $|\xi|^2$ near the origin, it
follows from Corollary 37.6 of \cite{Sa} that $Y$ is recurrent when
$d=2$ and transient when $d\ge 3$. By using the smoothness of the
density, one can easily check that, when $d\ge 3$, the Green
function of $Y$
$$
G^Y(x, y)=\int^{\infty}_0p^Y(t, x, y)dt
$$
is continuous off the set $\{(x, x): x\in \R^d\}$.

For any open set $D$, we use $\tau^X_D$ to denote the first exit
time from $D$ by the process $X$ and use $X^D$ to denote the process
obtained by killing the symmetric $\alpha$-stable process upon
leaving $D$. The process $X^D$ is usually called the killed
symmetric $\alpha$-stable process in $D$. The Dirichlet form of
$X^D$ is $({\cal E}, {\cal F}^D)$, where
$$
{\cal F}^D=\{u\in {\cal F}: u=0 \mbox{ on $D^c$ except for
a set of zero capacity }\}.
$$
For any $u, v\in {\cal F}^D$,
$$
{\cal E}(u, v)=\frac12\int_D\int_D(u(x)-u(y))(v(x)-v(y))J(x, y)dxdy
+\int_Du(x)v(x)\kappa_D(x)dx,
$$
where
$$
J(x, y):={\cal A}(d, -\alpha)|x-y|^{-(d+\alpha)}
\quad
\mbox{and}
\quad
\kappa_D(x):={\cal A}(d, -\alpha)\int_{D^c}|x-y|^{-(d+\alpha)}dy.
$$

Similarly, for any open set $D$, we use $\tau^Y_D$ to denote the
first exit time from $D$ by the process $Y$ and  use $Y^D$ to denote
the process obtained by killing the process $Y$ upon exiting $D$.
The Dirichlet form of $Y^D$ is $({\cal Q}, {\cal F}^D)$. For any $u,
v\in {\cal F}^D$,
$$
{\cal Q}(u, v)=\frac12\int_D\int_D(u(x)-u(y))(v(x)-v(y))J^Y(x, y)dxdy
+\int_Du(x)v(x)\kappa^Y_D(x)dx,
$$
where
\begin{equation}\label{e:J}
J^Y(x, y):=J(x, y)1_{\{|x-y|< 1\}}\quad
\mbox{and}
\quad
\kappa^Y_D(x):=\int_{D^c}J^Y(x, y)dy.
\end{equation}
Note that
\begin{equation}\label{cal}
0\le \kappa_{D}(x)-\kappa^Y_{D}(x)=
\int_{D^c\cap\{|x-y|\ge 1\}}J(x, y)dy\le
\int_{\{|x-y|\ge 1\}}J(x, y)dy=:{\cal B}(d, \alpha),
\quad \forall x\in D.
\end{equation}

Using the continuity of $p^Y$, it is routine (see, for instance,
the proof of Theorem 2.4 \cite{CZ} ) to show that $Y^D$ has a continuous
and symmetric density $p^Y_D(t, x, y)$. From this one can easily
show that the Green function $G^Y_D$ of $Y^D$ is continuous on
$(D\times D)\setminus\{(x, x): x\in D\}$.

\section{Comparability between Green Functions of Stable and
Truncated Stable Processes in small sets}

In this section we
take an open set $D$  with $\mbox{diam}(D) \le \frac12$, where
$\mbox{diam}(D)$  stands for the diameter of $D$.
In this case, we have from (\ref{e:J}) that
$J(x, y)=J^Y(x, y)$  for all $x, y\in D$.
So we can regard $Y^{D}$ as a Feynman-Kac
transform of the process $X^{D}$
with the potential $q_D(x):=\kappa_{D}(x)-\kappa^Y_{D}(x)$, that is,
the Feynman-Kac semigroup $(Q^D_t)$ defined by
$$
Q^{D}_tf(x)= \E_x\left[\exp(\int^t_0q_D(X^{D}_s)ds)f(X^{D}_t)\right]
$$
is the semigroup of $Y^D$.
Recall that $G^Y_{D}$ is the Green function of $Y^D$.
Let $G_{D}$ be the Green function of
$X^{D}$.
Since $q_D$ is nonnegative,
we see that
$G_{D}(x, y)\le G^Y_{D}(x, y)$ for all $x, y \in D$.
To get an upper bound for $G^Y_D$, we need to assume that
$D$ is $\kappa$-fat. We first recall the definition of
$\kappa$-fat set from \cite{SW}.

\begin{defn}\label{fat}
Let $\kappa \in (0,1/2]$. We say that an open set $D$ in
$\R^d$ is $\kappa$-fat if there exists $R>0$ such that for each
$Q \in \partial D$ and $r \in (0, R)$,
$D \cap B(Q,r)$ contains a ball $B(A_r(Q),\kappa r)$.
The pair $(R, \kappa)$ is called the characteristics of
the $\kappa$-fat open set $D$.
\end{defn}

Note that all Lipschitz domain and all non-tangentially accessible
domain (see \cite{JK} for the definition) are $\kappa$-fat.
Moreover, every {\it John domain} is  $\kappa$-fat (see Lemma 6.3 in
\cite{MV}). The boundary of a $\kappa$-fat open set can be highly
nonrectifiable and, in general, no regularity of its boundary can be
inferred.  Bounded $\kappa$-fat open set may be disconnected.

Suppose further that $D$ is a $\kappa$-fat set.
Since $q_D$ is bounded, one can use the 3G-type estimates for symmetric
stable processes in \cite{SW} to check that $q_D$ belongs to the
class ${\bf S}_{\infty}(X^D)$ there, thus
it follows from \cite{CS2} or \cite{CS3}
that $G^Y_{D}(x, y)$ is continuous and there exists a positive
constant $C$ such that
\begin{equation}\label{e:G_0}
G_{D}(x, y)\,\le\,G^Y_{D}(x, y)\,\le\, CG_{D}(x, y),
\quad \forall x, y \in D.
\end{equation}
But the constant $C$ above might depend on $D$. For later
applications, we will need the constant $C$ to be invariant under
scaling and translation. First we consider the case of  balls.

\begin{prop}\label{G_1}
There exists a positive constant $r_0\le \frac14$ such that for all
$r\in (0, r_0]$ and $a \in \R^d$, we have
$$
G_{B(a, r)}(x, y)\,\le\, G^Y_{B(a, r)}(x, y)\,\le\, 2\, G_{B(a,
r)}(x, y), \quad x, y\in B(a, r).
$$
\end{prop}

\pf Let $B_r:=B(0,r)$ with $r \le \frac14$. For any $z \in B_r$, let
$(\P^z_x , X^{B_r}_t)$ be the $G_{B_r}(\cdot , z)$-transform of
$(\P_x , X^{B_r}_t)$, that is, for any nonnegative Borel functions
$f$ in $B_r$,
$$\E^z_x \left[f(X^{B_r}_t)\right] = \E_x \left[\,
\frac{G_{B_r}(X^{B_r}_t , z)}
{G_{B_r}(x,z)}f(X^{B_r}_t)\right].
$$
It is well known that there exists a positive constant $C$
independent of $r$ such that
\begin{equation}\label{3g4balls}
\frac{G_{B_r}(x, y)G_{B_r}(y, z)}
{G_{B_r}(x, z)}\,\le\, C\,(|x-y|^{\alpha-d}+|y-z|^{\alpha-d}),
\quad \forall\, x, y, z\in B_r.
\end{equation}
So there exists a positive constant $r_0$ such that for any $r\in (0, r_0]$
and all $x, z\in B_r$,
\begin{equation}\label{e:3G_ball}
{\cal B}(d, \alpha)\,\E^z_x\tau^X_{B_r}\, = \,{\cal B}(d,
\alpha)\int_{B_r} \frac{G_{B_r}(x, y)G_{B_r}(y, z)} {G_{B_r}(x,
z)}dy \,<\frac12,
\end{equation}
 where ${\cal B}(d, \alpha)$ is the constant in (\ref{cal}).
Hence by (\ref{cal}) and Khasminskii's lemma (see, for instance,
Lemma 3.7 in \cite{CZ})
we get that for  $r \in (0, r_0 \wedge \frac14]$
$$
\E^z_x\left[\exp(\int^{\tau^X_{B_r}}_0q(X^{B_r}_s)ds)\right]\,\le\,
\E^z_x  \left[ \exp\left({\cal B}(d, \alpha)
\tau^X_{B_r}\right)\right]\,\le\, 2.
$$
Since
$$
G^Y_{B_r}(x, z)\,=\,G_{B_r}(x, z)\,\E^z_x\left[
\exp(\int^{\tau^X_{B_r}}_0q(X^{B_r}_s)ds)\right], \quad x, z\in B(0,
r),
$$
using the translation invariance property of our Green functions, we
arrive at our desired result. \qed

The 3G-type estimate for symmetric stable processes on $\kappa$-fat open sets
was proved in \cite{SW}.
It is easy to see that the constant appearing in the 3G estimate depends only
on  the characteristics  of the $\kappa$-fat open set
and the diameter of the set.
Moreover, by the scaling and translation invariant property of $X$,
the constant is invariant under scaling and translation of $D$

\begin{thm}\label{3G} (Theorem 6.1 in \cite{SW})
For a bounded $\kappa$-fat open set $D$ in $\R^d$,
there exists a constant $c$  depending only on the
characteristics $(\kappa, R)$ of $D$ and ${\rm diam}(D)$
such that
for $x, y, z \in D_r^a:=a+rD$,
\begin{equation}\label{3G_est}
 \frac{ G_{D_r^a}(x,y) G_{D_r^a}(y,z)}
{ G_{D_r^a}(x,z)} \,\le\, c\,\frac{ |x-z|^{d-\alpha}}
{|x-y|^{d-\alpha}  |y-z|^{d-\alpha}}.
\end{equation}
\end{thm}

By using \eqref{3G_est} instead of \eqref{3g4balls} in the proof of
Proposition \ref{G_1}, we immediately get the following result.

\begin{prop}\label{G_2}
Assume that $D$ is a bounded $\kappa$-fat open set in $\R^d$ with
 the characteristics  $(\kappa,R)$.
Then there exists constant  $r_1
=r_1(\kappa,R,\alpha,d,\mbox{diam}(D))  \le \frac12$diam$(D)$
 such that
for all $r\in (0, r_1]$ and $a \in \R^d$, we have
$$
G_{D_r^a}(x, y)\,\le\, G^Y_{D_r^a}(x, y)\,\le\, 2\,G_{D_r^a}(x, y),
\quad x, y\in D_r^a:=a+rD.
$$
\end{prop}

\pf We omit the details. \qed

\section{Harnack Inequality for Truncated Stable Processes}

In this section we will prove a Harnack inequality for
truncated stable processes.
It is well-known (see Lemma 6 of \cite{B})
that for any bounded Lipschitz domain $D$ in $\R^d$
(see Section 5 for the definition),
\begin{equation}\label{e:par}
\P_x(X_{\tau_D} \in \partial D)=0, \quad x \in D.
\end{equation}
The process $X$ has a L\'evy system $(N, H)$ with
$N(x, dy)={\cal A}(d, -\alpha)|x-y|^{-(d+\alpha)}dy$
and $H_t=t$
(see \cite{FOT}). Using this and (\ref{e:par})
we know that for every bounded Lipschitz domain $D$ and $f \ge 0$, we have
\begin{equation}\label{s_Levy}
\E_x\left[f(X_{\tau_D}) \right]
=\int_{\overline{D}^c}K_D(x, z)f(z)dz,\quad x \in D
\end{equation}
where
\begin{equation}\label{s2}
K_D(x, z)= {\cal A}(d,
-\alpha)\int_{D} \frac{G_D(x,y)}{|y-z|^{d+\alpha}} dy.
\end{equation}
Recall that for any ball $B(x, r)$, we use
$\tau^X_{B(x, r)}$ and $\tau^Y_{B(x, r)}$ to denote the first exit
time from the ball $B(x, r)$ by the processes $X$ and $Y$
respectively. Using Proposition \ref{G_1}, we can easily see that,
for $r\le r_0$,
$$
\E_0\tau^X_{B(x, r)}\,\le\, \E_0\tau^Y_{B(x, r)}\,\le\,
2\,\E_0\tau^X_{B(x, r)}.
$$
Thus it follows from Theorem 1 of \cite{Sz1} that for any
bounded Lipschitz domain $D$ in $\R^d$ we have
\begin{equation}\label{e:pary}
\P_x(Y_{\tau_D} \in \partial D)=0, \quad x \in D.
\end{equation}
The process $Y$ has a L\'evy system $(N^Y, H^Y)$ with
$N^Y(x, dy)={\cal A}(d, -\alpha)|x-y|^{-(d+\alpha)}1_{|x-y|<1}dy$
and $H^Y_t=t$
(see \cite{FOT}).
Using this and (\ref{e:pary}) we have the following result.

\begin{thm}\label{Poisson}
Suppose that $D$ is a bounded Lipschitz domain in $\R^d$.
Then there is a Poisson kernel
$K_D^Y(x,z)$ defined on $D \times \overline{D}^c_1$
such that
$$
\E_x\left[f(Y_{\tau_D})\right]
\,=\,
\int_{\overline{D}^c_1} K_D^Y(x,z) f(z)dz, \quad x \in D
$$
for every $f \ge 0$ on $\overline{D}^c_1$, where
$\overline{D}^c_1:=\{y \in \overline{D}^c : \mbox{\rm dist}(y, D) <
1\}$. Moreover,
$$
 K_D^Y(x,z)\,=\,
{\cal A}(d, -\alpha)\int_{D \cap \{ |y-z| <1\}}
 \frac{G_D^Y(x,y)}{|y-z|^{d+\alpha}}dy,
\quad (x,z) \in D \times \overline{D}^c_1.
$$
\end{thm}

\medskip

Using the L\'{e}vy system for $Y$ again, we know that for every
bounded open subset $D$ and  every $f \ge 0$ and $x \in D$,
\begin{equation}\label{newls}
\E_x\left[f(Y_{\tau_D});\,Y_{\tau_D-} \not= Y_{\tau_D}  \right]
=\int_{\overline{D}^c_1}
{\cal A}(d, -\alpha)\int_{D \cap \{ |y-z| <1\}}
\frac{G_D^Y(x,y)}{|y-z|^{d+\alpha}}dy f(z)dz.
\end{equation}
For notational convenience, we define
\begin{equation}\label{PK}
K_D^Y(x,z)\,:=\,
{\cal A}(d, -\alpha)\int_{D \cap \{ |y-z| <1\}}
\frac{G_D^Y(x,y)}{|y-z|^{d+\alpha}} dy,
\quad (x,z) \in D \times \overline{D}^c_1
\end{equation}
even if $D$ is not a bounded Lipschitz domain in $\R^d$, so
\eqref{newls} can be simply written as
$$
\E_x\left[f(Y_{\tau_D});\,Y_{\tau_D-} \not= Y_{\tau_D}  \right]
=\int_{\overline{D}^c_1} K_D^Y(x,z)f(z)dz.
$$

Recall that $r_0$ is the constant from Proposition \ref{G_1} and
$K_{B(x_0,r)}(x,z)$ is the
Poisson kernel of $B(x_0,r)$ with respect to $X$.
Let
$A(x, r, R):=\{ y \in \R^d: r \le |y-x| <R \}.$

\begin{lemma}\label{l:P_0}
Suppose that $x_0 \in \R^d$.
Then for every $r < \frac14$ and $z \in A(x_0, r, 1-r)$,
\begin{equation}\label{CC1}
K_{B(x_0,r)}(x,z) \,\le\, K_{B(x_0,r)}^Y(x,z),\qquad x \in B(x_0,r).
\end{equation}
If $r < r_0$ and $z \in A(x_0, r, \infty)$, then
\begin{equation}\label{CC2}
K_{B(x_0,r)}^Y(x,z) \,\le\, 2 K_{B(x_0,r)}(x,z),\qquad x \in B(x_0,r).
\end{equation}
\end{lemma}

\pf
Note that if $z \in A(x_0, r, 1-r)$ and $y \in B(x_0,r)$, then
$ |y-z| \le |x_0-y| +|x_0-z| < r + 1-r =1$.
Thus by Theorem \ref{Poisson},
$$
K_{B(x_0,r)}^Y(x,z)\,=\,
{\cal A}(d, -\alpha)\int_{B(x_0,r)}
 \frac{G_{B(x_0,r)}^Y(x,y)}{|y-z|^{d+\alpha}}
 dy.
$$
So (\ref{CC1}) follows from (\ref{e:G_0}) and (\ref{s_Levy}).

On the other hand, if $r < r_0$, by Theorem \ref{Poisson}
$$
K_{B(x_0,r)}^Y(x,z)\,=\,
{\cal A}(d, -\alpha)\int_{B(x_0,r) \cap \{ |y-z| <1\}}
 \frac{G_{B(x_0,r)}^Y(x,y)}{|y-z|^{d+\alpha}}
 dy\,\le\,
{\cal A}(d, -\alpha)\int_{B(x_0,r)}
 \frac{G_{B(x_0,r)}^Y(x,y)}{|y-z|^{d+\alpha}}
 dy.
$$
for every $z \in A(x_0, r, \infty)$ and $x \in B(x_0,r)$. Thus
(\ref{CC2}) follows from Proposition \ref{G_1} and
(\ref{s_Levy}).\qed

It is well-known that
\begin{equation}\label{P_f}
 K_{B(x_0,r)}(x,z)\,=\,c_1\,
\frac{(r^2-|x-x_0|^2)^{\frac{\alpha}2}}{(|z-x_0|^2-r^2)^{\frac{\alpha}2}}
\frac1{|x-z|^d}
\end{equation}
for some constant $c_1=c_1(d, \alpha) > 0$.

\begin{lemma}\label{l:P_1}
Suppose that $r < r_0$.
Then there exists a constant $c=c(d,\alpha) >0 $
such that  for any $z \in A(x_0, r, 1-r)$
and $x_1, x_2 \in B(x_0,\frac{r}{2})$,
$$
 c^{-1}K_{B(x_0,r)}(x_2,z) \,\le\, K_{B(x_0,r)}^Y(x_1,z)
\,\le\, c K_{B(x_0,r)}(x_2,z).
$$
\end{lemma}

\pf
By the previous lemma, we have for every $z \in A(x_0, r, 1-r)$,
\begin{equation}\label{aaa1}
K_{B(x_0,r)}(x,z) \,\le\, K_{B(x_0,r)}^Y(x,z) \,\le
\,2\, K_{B(x_0,r)}(x,z), \quad x \in B(x_0,r).
\end{equation}
By the explicit formula for  $K_{B(x_0,r)}(x,z)$ in (\ref{P_f}), we
see that there exist a constant $c_1=c_1(d,\alpha)$ such that for $x
\in B(x_0,\frac{r}{2})$ and $z \in A(x_0, r, 1-r)$
\begin{equation}\label{aaa2}
c_1^{-1} K_{B(x_0,r)}(x,z) \le  K_{B(x_0,r)}(x_0,z) \le c_1
K_{B(x_0,r)}(x,z).
\end{equation}
The inequalities (\ref{aaa1})-(\ref{aaa2}) imply the lemma.
\qed

\begin{lemma}\label{l:P_upper}
Suppose that $r < r_0$.
Then there exists a constant
$c=c(\alpha,d) > 0$ such that for any $z \in A(x_0, 1-r, 1+ r)$
and $x \in B(x_0,r)$
we have $K_{B(x_0,r)}^Y(x,z) \,\le\, c r^\alpha$.
\end{lemma}

\pf Without loss of generality, we may assume $x_0=0$. Fix  $r <
r_0$,  $z \in A(x_0, 1-r, 1+r)$ and  $x \in B(x_0,r)$. By (\ref{s2})
and (\ref{CC2}), we have
$$
K_{B(0,r)}^Y(x,z) \,\le\, c_1\int_{B(0,r)}
\frac{G_{B(0,r)}(x,y)}{|y-z|^{d+\alpha}}dy
$$
 for some constant $c_1=c_1(d, \alpha)$. Note that for any $y\in B(0,
r)$, $|y-z|\,\ge\,  |z| -|y|\,\ge\,1-r-r \,\ge\, \frac12.$ Since
$G_{B(0,r)}(x,y) \le c_2 |x-y|^{-d+\alpha}$ for some constant
$c_2=c_2(d, \alpha)$ (see, for instance, \cite{CS1}),
 we have
$$
K_{B(0,r)}^Y(x,z) \,\le \,c_1c_2\,\int_{B(0,r)}
\frac{dy}{|x-y|^{d-\alpha}} \,\le\,c_1c_2\, \int_{B(0,2r)}
\frac{dw}{|w|^{d-\alpha}} \le c_3 \,r^\alpha
$$
for some constant $c_3=c_3(d, \alpha)$. \qed

\begin{lemma}\label{l:P_2}
Suppose that $r < r_0$.
Then there exists
$c=c(\alpha,d) >0 $ such that for any $z \in A(x_0, 1-r, 1+\frac12 r)$
and $x \in B(x_0,\frac{r}{4})$,
$c^{-1} r^\alpha \le K_{B(x_0,r)}^Y(x,z) \le c r^\alpha.$
\end{lemma}

\pf
Without loss of generality, we may assume $x_0=0$.
Fix  $r < r_0$,  $z \in A(0, 1-r, 1+\frac12 r)$ and $x \in B(0,\frac{r}{4})$
. Let $B_r:=B(0,r)$.
By Lemma \ref{l:P_upper},  we only need to prove the lower bound.
Note that for any $y\in B_r$,
$
|y-z| \,\le\, |y|+|z|
<1 +\frac12 r +r \,<\,2.
$
So by
Theorem \ref{Poisson} and (\ref{e:G_0}), we have
$$
K_{B_r}^Y(x,z) \ge
c_1\int_{B_r \cap \{ |y-z| <1   \}}
G_{B_r}(x,y)dy \ge
c_1\int_{ \{  |y-z| <1,  |y| < \frac{7r}{8}  \}}
G_{B_r}(x,y)dy
$$
for some constant $c_1=c_1(d, \alpha)$.
Using the Green function estimates in \cite{CS1},
there exists $c_2=c_2(d, \alpha)$ such that
$$
 G_{B_r}(x,y) \ge c_2 |x-y|^{-d+\alpha}, \quad y \in B(0,\frac{7r}{8}).
$$
So we have
$$
K_{B_r}^Y(x,z) \ge c_1c_2
\int_{ \{ |y-z| <1 , |y| < \frac{7r}{8}  \}}
|x-y|^{-d+\alpha} dy.
$$
In the above integral, we will consider
the smallest possible open set to integrate on.
Let $z_r:=(0, \cdots, 0, 1+ \frac{r}{2}).$
The above integral is larger than or equal to
$$
\int_{ \{ |y-z_r| <1 , |y| < \frac{7r}{8}  \}}
|x-y|^{-d+\alpha} dy.
$$
Since
$|y| \ge  |z_r|-|y-z_r| > (1+ \frac{r}{2})-1= \frac{r}{2}$ for $|y-z_r| <1$,
we have
\begin{equation}\label{e:aa1}
|x-y| \le |x|+|y| < \frac14 r + |y| < 2|y|
\end{equation}
and
\begin{equation}\label{e:aa2}
\left\{ |y-z_r| <1 ,\, |y| < \frac{7r}{8}  \right\}
\, = \,\left\{ |y-z_r| <1 ,\, \frac{r}{2}< |y| < \frac{7r}{8}  \right\}.
\end{equation}
By a direct computation, one can show that
\begin{equation}\label{e:aa3}
\left\{ \frac{9r}{16}<y_d<\frac{11r}{16}, \,|(y_1, \cdots, y_{d-1})|
< \frac{r}{2}  \right\}
\, \subset\, \left\{ |y-z_r| <1 ,\, \frac{r}{2}< |y| < \frac{7r}{8}  \right\}.
\end{equation}
Putting (\ref{e:aa1}) and (\ref{e:aa3}) together, we get
$$
K_{B_r}^Y(x,z) \,\ge\, c_1c_2\,
\int_{\frac{9r}{16}}^{\frac{11r}{16}} \int_{\{|(y_1, \cdots,
y_{d-1})| < \frac{r}{2}  \}  }
|y|^{-d+\alpha}
 dy_1\cdots dy_d \,\ge\, c_3 \,r^\alpha
$$
for some constant $c_3=c_3(d, \alpha)$.
We have proved the lower bound.
\qed

Combining Lemmas \ref{l:P_1}-\ref{l:P_2}, we have proved the following.

\begin{lemma}\label{l:P_6_1}
Suppose that $r < r_0$.
Then there exists a constant $c=c(d,\alpha) >0$
such that  for any $z \in A(x_0, r, 1+\frac{r}{2})$
and $x_1, x_2 \in B(x_0,\frac{r}{4})$,
$$
 c^{-1}\,K_{B(x_0,r)}(x_2,z) \,\le\, K_{B(x_0,r)}^Y(x_1,z)\,
\le\, c\, K_{B(x_0,r)}(x_2,z).
$$
\end{lemma}

\begin{defn}\label{def:har1}
Let $D$ be an open subset of $\R^d$.
A function $u$ defined on $\R^d$ is said to be

\begin{description}
\item{(1)}  harmonic in $D$ with respect to $Y$ if
$$
\E_x\left[|u(Y_{\tau_{B}})|\right] <\infty
\quad \hbox{ and } \quad
u(x)= \E_x\left[u(Y_{\tau_{B}})\right],
\qquad x\in B,
$$
for every open set $B$ whose closure is a compact
subset of $D$;

\item{(2})
regular harmonic in $D$ with respect to $Y$ if it is harmonic in $D$
with respect to $Y$ and
for each $x \in D$,
$$
u(x)= \E_x\left[u(Y_{\tau_{D}})\right].
$$
\end{description}
\end{defn}

We define (regular) harmonic function with respect to $X$
similarly. The next lemma is a preliminary version of the
Harnack inequality for $Y$
and it is an immediate consequence of Lemma \ref{l:P_1}.

\begin{lemma}\label{T:Har_1}
Suppose that $r\le r_0$. There exists a constant $c=
c(d, \alpha)$ such that
$$
c^{-1}u(y)\le u(x)\le cu(y), \quad y\in B(x, \frac{r}2)
$$
for any nonnegative function $u$ which is regular
harmonic in $B(x, r)$ and zero in $B(x, 2r)^c$.
\end{lemma}

Now we are ready to prove a (scale-invariant) Harnack inequality for $Y$.

\begin{thm}\label{T:Har}
Suppose $x_{1}, x_{2}\in \R^d$, $r < r_0$ are
such that $|x_{1}-x_{2}|< Mr$
for some $M \le \frac{1}{r} - \frac12$.
Then there exists a constant $J >0 $ depending only
on $d$ and $\alpha$, such that
$$
J^{-1}M^{-(d+\alpha)}u(x_{2})\,\leq \,u(x_{1})
\,\leq\, JM^{d+\alpha}u(x_{2})\,
$$
for every nonnegative function $u$
which is regular harmonic with respect to $Y$ in
$B(x_{1}, r)\cup B(x_{2},r)$.
\end{thm}

\pf Fix  $r <r_0$,  $x_{1}, x_{2}\in \R^d$ and
a nonnegative regular harmonic function $u$ in
$B(x_{1}, r)\cup B(x_{2},r)$ with respect to $Y$.
Let $B^i=B(x_i, \frac14 r)$, $i=1,2$.

We split into two cases.
First we deal with the case $ |x_{1}- x_{2}|  <  \frac14 r$.
In this case we have $ \emptyset \neq B^1\cap B^2 \supset \{x_1,x_2 \}$.
By Theorem \ref{Poisson} we have for any $y\in B^1\cap B^2$ and $i=1, 2$,
\begin{eqnarray*}
&&u(y) = \E_{y} \left[u(Y_{\tau_{B(x_i, r)}})\right]\\
&&= \int_{A(x_i, r, 1+\frac{r}{2})} K_{B(x_i,r)}^Y(y,z)u(z) dz +
\E_{y} \left[u(Y_{\tau_{B(x_i, r)}});\, Y_{\tau_{B(x_i, r)}}
\in   A(x_i, 1+\frac{r}{2}, 1+r)    \right].
\end{eqnarray*}
By Lemma \ref{l:P_6_1},
\begin{eqnarray*}
 \int_{A(x_i, r, 1+\frac{r}{2})} K_{B(x_i,r)}^Y(y,z) u(z)dz
&\le& c_1 \int_{A(x_i, r, 1+\frac{r}{2})} K_{B(x_i,r)}^Y(x_i,z) u(z)dz\\
&=&c_1\, \E_{x_i} \left[u(Y_{\tau_{B(x_i, r)}});\, Y_{\tau_{B(x_i, r)}}
\in   A(x_i, r, 1+\frac{r}{2})    \right],
\end{eqnarray*}
for some constant $c_1=c_1(d, \alpha)$.
Note that by Theorem \ref{Poisson}
\begin{eqnarray*}
&&\P_{y} \left( Y_{\tau_{B(x_i, r)}}
\in   A(x_i, 1+\frac{r}{2}, 1+r) ,\, \tau_{B(x_i, r)} =
\tau_{B(x_i, \frac{r}{2})}\right) \\
&& = \P_{y} \left( Y_{\tau_{B(x_i, \frac{r}{2})}}
\in   A(x_i, 1+\frac{r}{2}, 1+r)\right)
\, =\,  \int_{A(x_i, 1+\frac{r}{2}, 1+r )} K_{B(x_i,\frac{r}{2})}^Y(y,z) dz=0.
\end{eqnarray*}
Thus by the strong Markov property, we have
\begin{eqnarray*}
&&\E_{y} \left[\,u(Y_{\tau_{B(x_i, r)}});\, Y_{\tau_{B(x_i, r)}}
\in   A(x_i, 1+\frac{r}{2}, 1+r)     \right]\\
&&=\E_{y} \left[\,u(Y_{\tau_{B(x_i, r)}});\, Y_{\tau_{B(x_i, r)}}
\in   A(x_i, 1+\frac{r}{2}, 1+r) ,\,    \tau_{B(x_i, r)} >
\tau_{B(x_i, \frac{r}{2})}   \right]\\
&&=
\E_{y} \left[\,    \E_{Y_{\tau_{B(x_i, \frac{r}{2})}}} \left[
u(Y_{\tau_{B(x_i, r)}});\, Y_{\tau_{B(x_i, r)}}
\in   A(x_i, 1+\frac{r}{2}, 1+r)   \right]
 1_{ A(x_i, \frac{r}{2}, r)} (Y_{\tau_{B(x_i, \frac{r}{2})}})
 \right].
\end{eqnarray*}
For $i=1, 2$, let
$$
g_i(z):=\E_{z} \left[    u(Y_{\tau_{B(x_i, r)}});\, Y_{\tau_{B(x_i, r)}}
\in   A(x_i, 1+\frac{r}{2}, 1+r)   \right]
$$
for $z \in  A(x_i, \frac{r}{2}, r)$, and zero otherwise.
Then we have from the above argument that
$$
\E_{y} \left[u(Y_{\tau_{B(x_i, r)}});\, Y_{\tau_{B(x_i, r)}}
\in   A(x_i, 1+\frac{r}{2}, 1+r)     \right]=
\E_{y} \left[g_i(Y_{\tau_{B(x_i, \frac{r}{2})}})
 \right].
$$
Since, for $i=1, 2$, the function
$y\mapsto \E_{y} \left[g_i(Y_{\tau_{B(x_i, \frac{r}{2})}})\right]$
is regular harmonic on $B(x_i, \frac{r}{2})$ with respect to $Y$,
and is zero on $\overline{B(x_i, r)}^c$, we get
by Lemma \ref{T:Har_1} that for $y\in B^1\cap B^2$,
\begin{eqnarray*}
&&\E_{y} \left[u(Y_{\tau_{B(x_i, r)}});\, Y_{\tau_{B(x_i, r)}}
\in   A(x_i, 1+\frac{r}{2}, 1+r)     \right]\\
&&\le c_2\,\E_{x_i} \left[    \E_{Y_{\tau_{B(x_i, \frac{r}{2})}}} \left[
u(Y_{\tau_{B(x_i, r)}});\,
Y_{\tau_{B(x_i, r)}}
\in   A(x_i, 1+\frac{r}{2}, 1+r)   \right]\,
 1_{ A(x_i, \frac{r}{2}, r)} (Y_{\tau_{B(x_i, \frac{r}{2})}})
 \right]\\
&&=\,c_2\,\E_{x_i} \left[u(Y_{\tau_{B(x_i, r)}});\, Y_{\tau_{B(x_i, r)}}
\in   A(x_i, 1+\frac{r}{2}, 1+r)     \right],
\end{eqnarray*}
for some constant $c_2=c_2(d, \alpha)$.
Combining the two parts together, we get
that
\begin{equation}\label{e:H11}
u(y)
\,\le\, c_3
\,u(x_i), \qquad y\in B^1\cap B^2
\end{equation}
for some constant $c_3=c_3(d, \alpha).$
Therefore
$c_3^{-1} u(x_{2})\leq u(x_{1})
\leq c_3 u(x_{2})$.

Now we consider the case when $\frac14 r \le |x_{1}- x_{2}| \le M r$
with $M \le \frac{1}{r} - \frac12$.
Since $M \le \frac{1}{r} - \frac12$,
we have
$|x_{1}-x_{2}|< 1-\frac12 r$, and
$|y-w| \le |y-x_2|+|w-x_1|+ |x_{1}-x_{2}| <1$
for $(y,w) \in B(x_2, \frac{r}8)\times B(x_1, \frac{r}8)$.
Thus, by Proposition \ref{G_1} and Theorem \ref{Poisson}, we have for
$w \in B(x_1, \frac{r}8)$,
\begin{eqnarray*}
K_{B(x_2, \frac{r}8)}^Y(x_2,w) &=&
{\cal A}(d, -\alpha)\int_{B(x_2, \frac{r}8)} \frac{G^Y_{B(x_2, \frac{r}8)}
(x_2,y)}
{|y-w|^{d+\alpha}} dy\\
&\ge &
{\cal A}(d, -\alpha)\int_{B(x_2, \frac{r}8)} \frac{G_{B(x_2, \frac{r}8)}
(x_2,y)}
{|y-w|^{d+\alpha}} dy
~=~
K_{B(x_2, \frac{r}8)}(x_2,w).
\end{eqnarray*}
From (\ref{P_f}), we have for
$w \in B(x_1, \frac{r}8)$,
$$
K_{B(x_2, \frac{r}8)}^Y(x_2,w) \,\ge\, c_4 \,
\frac{r^\alpha}{|x_2-w|^{d+\alpha}}
\,\ge\, c_4\, \frac{r^{-d}}{(2M)^{d+\alpha}}
$$
for some constant $c_4=c_4(d, \alpha)$,
because $|x_2-w| \le |x_1-x_2|+|w-x_1| < (M+\frac18)r \le 2Mr$.
For any $y\in B(x_1, \frac{r}8)$, $u$ is regular harmonic in
$B(y, \frac{7r}8)\cup B(x_1, \frac{7r}8)$. Since
$|y-x_1|< \frac{r}8$, we can apply the conclusion
of the first case with $x_2=y$ and $r$ replaced by $\frac{7r}8$ to get that
$$
u(y)\ge c_5 u(x_1), \quad y\in B(x_1, \frac{r}8),
$$
for some constant $c_5=c_5(d, \alpha)$.
Therefore
\begin{eqnarray*}
u(x_2) &=& \E_{x_2}\left[u(Y_{\tau_{B(x_2,\frac{r}8) }})\right]
  \ge \E_{x_2}\left[u(Y_{\tau_{B(x_2,\frac{r}8)}});
Y_{\tau_{B(x_2,\frac{r}8)}} \in B(x_1,\frac{r}8) \right]\\
  &\ge &c_5 \,u(x_1)\, \P_{x_2}\left(Y_{\tau_{B(x_2,\frac{r}8)}}
\in B(x_1,\frac{r}8) \right)
 = c_5\, u(x_1) \int_{B(x_1,\frac{r}8)} K_{B(x_2,\frac{r}8)}^Y(x_2,w)dw \\
 &\ge  &c_4c_5\, u(x_1) \int_{B(x_1,\frac{r}8)}
\frac{r^{-d}}{(2M)^{d+\alpha}}dw
 \,\,\ge\,\,  c_6\, u(x_1) \,M^{-(d+\alpha)},
\end{eqnarray*}
for some constant $c_6=c_6(d, \alpha)$.
We have thus proved the right hand side inequality in the conclusion
of the theorem. The inequality on the left hand side can be
proved similarly.
\qed

The Harnack inequality above  is similar to the Harnack inequality
(Lemma 2) for symmetric stable processes in \cite{B}, the difference
is that we have to require that the two balls are not too far apart.
Because our process can only make jumps of size at most 1, one can
easily see that, without the assumption above, the Harnack
inequality fails.

As a consequence of the theorem above we immediately get
the following

\begin{corollary}\label{c:Har_1}
Suppose that $r\le r_0$. There exists a constant $c=
c(d, \alpha) >0$ such that
$$
c^{-1}\,u(y)\,\le\, u(x)\,\le\, c\,u(y), \qquad y\in B(x, \frac{r}2)
$$
for any nonnegative function $u$ which is
harmonic in $B(x, r)$.
\end{corollary}

Using this and a standard chain argument, we can get the following

\begin{corollary}\label{c:Har_2}
Suppose that $D$ is a  domain (i.\/e., a connected open set) in $\R^d$
and $K$ is a compact subset of $D$. There exists a constant $c=
c(D, K, \alpha)>0$ such that
$$
c^{-1}\,u(y)\,\le\, u(x)\,\le\, c\,u(y), \qquad x, y\in K
$$
for any nonnegative function $u$ which is
harmonic in $D$
\end{corollary}

\section{Boundary Harnack Principle for Truncated Stable Processes}

In this section we  will prove two versions of the boundary Harnack principle
for truncated stable processes.
Throughout this section, $r_0$ is the  constant in
Proposition \ref{G_1}.

We will use $\AA_\alpha$ to denote the $L_2$-generator of $Y$, and
$C^{\infty}_c(\RR^d)$ to denote
the space of continuous function with compact support. It is well-known that
$C^{\infty}_c(\RR^d)$ is in   the domain of  $\AA_\alpha$ and, for every
$\phi \in C^{\infty}_c(\RR^d)$,
\begin{equation}\label{generator}
 \AA_\alpha  \phi(x)=\AA(d, -\alpha)\int_{|y|<1}
\frac{\phi(x+y)-\phi(x)-(\nabla \phi(x)\cdot y)1_{B(0, \eps)}(y)}
{|y|^{d+\alpha}}dy,
\end{equation}
(see Section 4.1 in \cite{Sk}).

For any $\lambda>0$, let $G^{Y, \lambda}(x,y)$ be the
$\lambda$-Green function of $Y$.
We have
$(\AA_{\alpha}-\lambda)
G^{Y, \lambda}(x,y)=-\delta_x(y)$ in the weak sense.
For any bounded open subset $D$ of $\RR^d$, let
$G^{Y, \lambda}_D(x, y)$ be the
$\lambda$-Green function of $Y^D$.
Since
$
G^{Y, \lambda}_D(x,y)=G^{Y, \lambda}(x,y)
-\E_x[e^{-\lambda \tau_D}G^{Y, \lambda}(Y_{\tau_D},y)]
$,
we have, by the symmetry of $\AA_\alpha$,  for any $x\in D$ and any
nonnegative $\phi \in C^{\infty}_c(\RR^d)$,
\begin{eqnarray*}
&&\int_D   G^{Y, \lambda}_D(x,y) (\AA_{\alpha}-\lambda)  \phi(y)dy
=\int_{\R^d} G^{Y, \lambda}_D(x,y) (\AA_{\alpha}-\lambda)  \phi(y)dy\\
 &&=
\int_{\R^d} G^{Y, \lambda}(x,y)  (\AA_{\alpha}-\lambda)  \phi(y)dy-
 \int_{\RR^d} \E_x[e^{-\lambda \tau_D}G^{Y, \lambda}(Y_{\tau_D},y)]
(\AA_{\alpha}-\lambda)  \phi(y)dy\\
 &&=
\int_{\R^d} G^{Y, \lambda}(z,y) (\AA_{\alpha}-\lambda) \phi(y)dy\\
&&\,\,\,\,\,\,\,\,\,- \int^{\infty}_0e^{-\lambda t}\int_{D^c}
\int_{\RR^d} G^{Y, \lambda}(z,y) (\AA_{\alpha}-\lambda) \phi(y)dy
\P_x(Y_{\tau_D} \in dz, \tau_D\in dt)\\
 &&=-\phi(x)+\int^{\infty}_0 e^{-\lambda t} \int_{D^c}
\phi(z)\P_x(Y_{\tau_D} \in dz,
\tau_D\in dt)
\,=\,-\phi(x)+\E_x[e^{-\lambda \tau_D}\phi(Y_{\tau_D})].
\end{eqnarray*}
In particular, if $\phi(x)=0$ for $x\in D$, we have
\begin{equation}\label{lambda_har_gen}
\E_x\left[e^{-\lambda \tau_D} \phi(Y_{\tau_D})\right] = \int_D
G^{Y, \lambda}_D(x,y)
(\AA_{\alpha}-\lambda)  \phi(y)dy.
\end{equation}
Since $G^{Y, \lambda}_D(x,y)$ increases to
$G^Y_D(x, y)$ as $\lambda\downarrow 0$ and $(\AA_{\alpha} -\lambda)
\phi$ is bounded for small $\lambda$, by
letting $\lambda\downarrow 0$ in the equation above,
the dominated convergence theorem gives
\begin{equation}\label{har_gen}
\E_x\left[\phi(Y_{\tau_D})\right] = \int_D
G^Y_D(x,y)
\AA_{\alpha}  \phi(y)dy
\end{equation}
for any $x\in D$ satisfying $\phi(x)=0$.
Take a sequence of radial functions $\phi_m$ in $C^{\infty}_c(\RR^d)$
such that $0\le \phi_m\le 1$,
\[
\phi_m(y)=\left\{
\begin{array}{lll}
0, & |y|<1/2\\
1, & 1\le |y|\le m+1\\
0, & |y|>m+2,
\end{array}
\right.
\]
and that $\sum_{i, j}|\frac{\partial^2}{\partial y_i\partial y_j}
\phi_m|$ is uniformly bounded.
Define $\phi_{m, r}(y)=\phi_m(\frac{y}{r})$ so that
$0\le \phi_{m, r}\le 1$,
\[
\phi_{m, r}(y)=\left\{
\begin{array}{lll}
0, & |y|<r/2\\
1, & r\le |y|\le r(m+1)\\
0, & |y|>r(m+2),
\end{array}
\right.
\]
and
$$
\sup_{y\in \RR^d}
\sum_{i, j}\left|\frac{\partial^2}{\partial y_i\partial y_j}
\phi_{m, r}(y)\right| \,<\, c\, r^{-2}.
$$
We claim that there exists a constant $C>0$ such that
for all $r\in (0, 1)$,
\begin{equation}\label{e2.1}
\sup_{m \ge 1}  \sup_{y\in \RR^d} |\AA_{\alpha}\phi_{m,r}(y)|\,\le\, C
r^{-\alpha}.
\end{equation}
In fact, we have
\begin{eqnarray*}
&&|\AA_{\alpha}\phi_{m,r}(x)|
\le \AA(d, -\alpha)
\left|\int_{\{|y|<1\}}
\frac{\phi_{m,r}(x+y)-\phi_{m,r}(x)-(\nabla \phi_{m,r}(x)
\cdot y)1_{B(0, r)}(y)}
{|y|^{d+\alpha}}dy \right|\\
&&= \AA(d, -\alpha)\left(\left|\int_{\{|y|\le r\}}
\frac{\phi_{m,r}(x+y)-\phi_{m,r}(x)-(\nabla \phi_{m,r}(x)\cdot y)}
{|y|^{d+\alpha}}dy\right|+\int_{\{r<|y|< 1\}}\frac1{|y|^{d+\alpha}}dy
\right)\\
&&\le \AA(d, -\alpha)\left(\frac{c}{r^2}\int_{\{|y|\le r \}}
\frac{|y|^2}{|y|^{d+\alpha}}dy +\int_{\{r<|y|< 1\}}\frac1{|y|^{d+\alpha}}dy
\right)
\le C_1r^{-\alpha},
\end{eqnarray*}
for some constant $C_1=C_1(d, \alpha)>0$.
When $D \subset B(0,r)$ for some $r\in (0, 1)$,
we get, by combining (\ref{har_gen})
and  (\ref{e2.1}), that for any $x\in D\cap B(0, \frac{r}2)$,
\[
\P_x\left(Y_{\tau_D} \in B(0, r)^c\right)\,=\,
\lim_{m\to \infty}\P_x\left(Y_{\tau_D}
\in A(0, r, (m+1)r)\right)
\,\le\, C\,r^{-\alpha}\int_D  G_D^Y(x,y)dy.
\]
We have proved the following.

\begin{lemma}\label{l2.1}
Let $r\in (0, 1)$ and $D$ be an open subset with $D\subset B(0, r)$.
Then
\[
\P_x\left(Y_{\tau_D} \in B(0, r)^c\right)
\,\le\, C\,r^{-\alpha}\int_D G_D^Y(x,y)dy, \quad x \in D\cap B(0, \frac{r}2)
\]
for some constant $C=C(d, \alpha)>0$.
\end{lemma}

Recall that $r_0$ is  the constant  from  Proposition \ref{G_1}.

\begin{lemma}\label{l2.1_1}
Let $D$ be an open set
such that $B(A, \kappa r)\subset D\subset B(0, r)$
for some $r>0$ and $\kappa\in (0, 1)$.
If  $ r < r_0$,  then
$$
\P_x\left(Y_{\tau_{D\setminus B(A, \kappa r) }} \in B(A, \kappa r) \right)
\,\ge\, C\, r^{-\alpha}\, \kappa^{d} \, \int_D
 G_D^Y(x,y)dy, \quad x \in D \setminus B(A, \kappa r)
$$
for some constant $C=C(d, \alpha)>0$.
\end{lemma}

\pf Fix a point $x\in D\setminus B(A, \kappa r)$ and let
$B:=B(A, \frac{\kappa r}2)$.
Since $G_D^Y(x,\,\cdot\,)$ is harmonic in $D\setminus
\{x\}$ with respect to $Y$,
\[
G_D^Y(x,A)\,=\,\int_{D\cap B^c}K_B^Y(A, y)G_D^Y(x,y)dy
\,\ge\,\int_{D\cap B(A, \frac{3\kappa r}4)^c}K_B^Y(A,
y)G_D^Y(x,y)dy.
\]
Since $r < \frac14$, by (\ref{CC1}),  we have
\[
G_D^Y(x,A)
\,\ge\,\int_{D\cap B(A, \frac{3\kappa r}4)^c}K_B(A,
y)G_D^Y(x,y)dy.
\]
Since $\frac{3\kappa r}4\le |y-A|\le 2r$
for $y\in B(A, \frac{3\kappa r}4)^c\cap D$,
it follows from (\ref{P_f}) that
\begin{eqnarray*}
G_D^Y(x,A)&\ge& c_1\int_{D\cap B(A, \frac{3\kappa r}4)^c}
\frac{(\kappa r)^{\alpha}}{|y-A|^{d+\alpha}}G_D^Y(x,y)dy\\
&\ge& c_2\kappa^{\alpha}r^{-d}\int_{D \cap B(A, \frac{3\kappa
r}4)^c}G_D^Y(x,y) dy,
\end{eqnarray*}
for some constants $c_1=c_1(d, \alpha)$ and $c_2=c_2(d, \alpha)$.
Applying Theorem \ref{T:Har} we get
\[
\int_{B(A, \frac{3\kappa r}4)} G_D^Y(x,y)   dy\le
c_3 \int_{B(A, \frac{3\kappa r}4)}  G_D^Y(x,A)dy
\,\le\,c_4\,\kappa^d\,r^d\,G_D^Y(x,A),
\]
for some constants $c_3=c_3(d, \alpha)$ and $c_4=c_4(d, \alpha)$.
Combining these two estimates we get that
\[
\int_{D} G_D^Y(x,y)   dy    \,\le\, c_5\,\kappa^{-\alpha}\,r^{d}\,
G_D^Y(x,A)
\]
for some constant $c_5=c_5(d, \alpha)$.

Let $\Omega=D\setminus B(A, \frac{\kappa r}2)$.
Since diam$(D) <1 $, from (\ref{PK}) we have
 for $z\in B(A, \frac{\kappa r}4)$
\begin{equation}\label{e:KK0}
K^Y_{\Omega}(x, z) =
{\cal A}(d, -\alpha)\int_{\Omega}
\frac{G_{\Omega}^Y(x,y)}{|y-z|^{d+\alpha}} dy
\end{equation}
Note that for any $z\in B(A, \frac{\kappa r}4)$ and $y\in \Omega$,
$
2^{-1}|y-z|\le |y-A|\le 2|y-z|.
$
Thus
we get from (\ref{e:KK0}) that for $z\in B(A, \frac{\kappa r}4)$,
\begin{equation}\label{e:KK1}
2^{-d-\alpha}K^Y_{\Omega}(x, A)  \,\le \,K^Y_{\Omega}(x, z)
\,\le\,
2^{d+\alpha}K^Y_{\Omega}(x, A).
\end{equation}
Using the harmonicity of $G_D^Y(\cdot, A)$ in $D\setminus\{A\}$
with respect to $Y$,
we can split $G_D^Y(\cdot, A)$ into two parts:
\begin{eqnarray*}
&&G^Y_D(x, A)
=\E_y \left[G^Y_D(Y_{\tau_{\Omega}},A)\right]\\
&&=\E_y \left[G^Y_D(Y_{\tau_{\Omega}},A):\,Y_{\tau_{\Omega}}
\in B(A, \frac{\kappa r}4)  \right]\,+\, \E_y
\left[G^Y_D(Y_{\tau_{\Omega}},A):\,Y_{\tau_{\Omega}}
\in \{\frac{\kappa r}4\le |y-A|\le \frac{\kappa r}2\}\right]\\
&& :=I_1+I_2.
\end{eqnarray*}
By the monotonicity of Green functions and Proposition \ref{G_1},
\begin{equation}\label{GGGG}
G_D^Y(y, A)\,\le\,G_{B(0,r)}^Y(y, A) \,\le\, 2\,G_{B(0,r)}(y, A)\,
\le\, 2\, G(y, A),
 \qquad y \in B(0,r),
\end{equation}
where $G(\cdot, \cdot)$ is the Green function of $X$.
So using (\ref{e:KK1}) twice and the explicit formula for
$G(\,\cdot\,, \,\cdot \,)$, we have
\begin{eqnarray*}
&&I_1 \,\le\, 2^{d+\alpha}\,K^Y_{\Omega}(x,A)
\int_{B(A, \frac{\kappa r}4)}G_D^Y(y, A)dy
\,\le\, c_6 \,K^Y_{\Omega}(x,A)
\int_{B(A, \frac{\kappa r}4)}\frac{dy}{|y-A|^{d-\alpha}}
 \\
&&\le\, c_7\kappa^{\alpha}r^{\alpha}K^Y_{\Omega}(x, A)
\,\le\, c_8\kappa^{\alpha-d}r^{\alpha-d}\int_{B(A, \frac{\kappa r}4)}
K_{\Omega}^Y(x, z)dz,
\end{eqnarray*}
for some constants $c_i=c_i(d, \alpha), i=6, 7, 8.$
On the other hand, by (\ref{GGGG})
\begin{eqnarray*}
&&I_2 \,\le\,
\int_{\{\frac{\kappa r}4\le |y-A|\le \frac{\kappa r}2\}}
G_{B(0,r)}^Y(y,A) \P_x(Y_{\tau_{\Omega}} \in dy)   \\
&&\le
 c_9\int_{\{\frac{\kappa r}4\le |y-A|\le \frac{\kappa r}2\}}
\frac1{|y-A|^{d-\alpha}}  \P_x(Y_{\tau_{\Omega}} \in dy)
\,\le\,c_{10}\,
\kappa^{\alpha-d}\,r^{\alpha-d} \, \P_x\left(Y_{\tau_{\Omega}} \in
\{\frac{\kappa r}4\le
|y-A|\le \frac{\kappa r}2\}\right),
\end{eqnarray*}
for some constants $c_i=c_i(d, \alpha), i=9, 10$.
Therefore
$$
G^Y_D(x, A)
\,\le\, c_{11}\,
\kappa^{\alpha-d}\,r^{\alpha-d}\,\P_x\left(Y_{\tau_{\Omega}} \in
B(A, \frac{\kappa r}2)\right).
$$
for some constant $c_{11}=c_{11}(d, \alpha)$.
This implies that
\[ \int_{D} G_D^Y(x,y)   dy
\,\le\,c_{12}\,\kappa^{-d}\, r^\alpha \,
\P_x\left(Y_{\tau_{D\setminus B(A, \frac{\kappa r}2)}}
\in B(A, \frac{\kappa r}2)\right),
\]
for some constant $c_{12}=c_{12}(d, \alpha)$. It follows immediately
that
$$
\int_{D} G_D^Y(x,y)   dy  \,
\le\,c_{12}\,\kappa^{-d} \,r^\alpha \,
\P_x\left(Y_{\tau_{D\setminus B(A, \kappa r)}}
\in B(A, \kappa r)\right).
$$
\qed

Combining Lemmas \ref{l2.1}-\ref{l2.1_1} and using
the translation invariant property, we have
the following

\begin{lemma}\label{l2.3}
Let $D$ be an open set
such that $B(A, \kappa r)\subset D\subset B(0, r)$
for some $r>0$ and $\kappa\in (0, 1)$.
If  $r < r_0$,   then
\[
\P_x\left(Y_{\tau_{D}} \in B(0, r)^c\right) \,\le\,
C\,\kappa^{-d}\,\P_x\left(Y_{\tau_{D\setminus B(A, \kappa r) }}
\in B(A, \kappa r) \right),
 \qquad x\in D\cap B(0,
\frac{r}2),
\]
for some constant $C=C(d, \alpha)>0$.
\end{lemma}

The next lemma is adapted from \cite{B} (see page 54-55 in \cite{B}).

\begin{lemma}\label{l2.U}
Let $D$ be an open set and $0<2r<r_0$.
For any positive function $u$,
there is a $\sigma\in (\frac{10}{6}r, \frac{11}{6}r)$
such that for any  $a \in (-2,\frac32]$, $z_0\in \R^d$ and $ x \in D
\cap B(z_0, \frac{3}{2}  r)$,
$$\E_x\left[u(Y_{\tau_{D \cap B(z_0, \sigma)}});
Y_{\tau_{D \cap B(z_0, \sigma)}} \in
A(z_0, \sigma, 1- a r) \right]
\le\, C\, r^{\alpha}
\int_{A(z_0, \frac{10r}6, 1- a r)}
\frac{u(y)}{|y|^{d+\alpha}}dy
$$
for some constant $C=C(d, \alpha)$.
\end{lemma}

\pf Without loss of generality, we may assume $z_0=0$. Note that
\begin{eqnarray*}
&&\int^{\frac{11}{6}r}_{\frac{10}{6}r}\int_{A(0, \sigma, 2r)}
(|y|-\sigma)^{-\frac{\alpha}2}
u(y)dyd\sigma
 =\int_{A(0, \frac{10}{6}r , 2r)}
\int^{ |y| \wedge \frac{11}{6}r}_{\frac{10}{6}r}
(|y|-\sigma)^{-\frac{\alpha}2}d\sigma u(y)dy \\
&&\le c_1
\int_{A(0, \frac{10r}6, 2r)} (|y|- \frac{10r}6)^{-\frac{\alpha}2+1} u(y)dy
\le c_1
\int_{A(0, \frac{10r}6, 2r)} |y|^{1-\frac{\alpha}2} u(y)dy,
\end{eqnarray*}
for some constant $c_1=c_1(\alpha)$.
Thus there is a
$\sigma\in (\frac{10}{6}r, \frac{11}{6}r)$ such that
\begin{equation}\label{e:int}
\int_{A(0, \sigma, 2r)}(|y|-\sigma)^{-\frac{\alpha}2}u(y)dy
\,\le\, c_2\,r^{-1}\int_{A(0, \frac{10r}6, 2r)}u(y) |y|^{1-\frac{\alpha}2} dy
\end{equation}
for some constant $c_2=c_2(\alpha)$.
Let  $x \in D \cap B(0,  \frac{3}{2}  r)$.
Note that, by Theorem \ref{Poisson}
\begin{eqnarray*}
&& \E_x\left[u(Y_{\tau_{D \cap B(0, \sigma)}}); Y_{\tau_{D \cap
B(0, \sigma)}} \in
A(0, \sigma, 1- a r) \right]\\
&&= \E_x\left[u(Y_{\tau_{D \cap B(0, \sigma)}}); Y_{\tau_{D \cap
B(0, \sigma)}} \in
A(0, \sigma, 1- a r), \, \tau_{D \cap
B(0, \sigma)} =\tau_{B(0, \sigma)}  \right]\\
&&= \E_x\left[u(Y_{\tau_{ B(0, \sigma)}}); Y_{ \tau_{B(0, \sigma)}} \in
A(0, \sigma, 1- a r), \, \tau_{D \cap B(0, \sigma)} =
\tau_{B(0, \sigma)}  \right]\\
&&\le \E_x\left[u(Y_{\tau_{ B(0, \sigma)}}); Y_{\tau_{B(0, \sigma)}} \in
A(0, \sigma, 1- a r)  \right]
\,=\,\int_{A(0, \sigma, 1- a r)}K^Y_{B(0, \sigma)}(x, y)u(y)dy.
\end{eqnarray*}
Since  $ \sigma  <2r < r_0$, by (\ref{CC2}) and (\ref{P_f})
we have
\begin{eqnarray*}
&& \E_x\left[u(Y_{\tau_{D \cap B(0, \sigma)}}); Y_{\tau_{D \cap
B(0, \sigma)}} \in
A(0, \sigma, 1- a r) \right]
\,\le\, 2
 \int_{  A(0, \sigma, 1- a r)   }
K_{B(0, \sigma)}(x, y)u(y)dy\\
&&\le\,
c_3 \left(\int_{A(0, \sigma, 2r)}+\int_{A(0, 2r , 1- a r  )}\right)
\frac{(\sigma^2-|x|^2)^{\frac{\alpha}2}}{(|y|^2-\sigma^2)^
{\frac{\alpha}2}}\frac1{|y-x|^d}u(y)dy
\end{eqnarray*}
for some constant $c_3=c_3(d, \alpha)$.
For $y \in A(0,  2r , 1- a r)   $,
$|y|^2-\sigma^2 \ge \frac1{12}|y|^2$
and $ \sigma^2-|x|^2 \le cr^2$. So
$$
\frac{(\sigma^2-|x|^2)^{\frac{\alpha}2}}{(|y|^2-\sigma^2)^
{\frac{\alpha}2}}\frac1{|y-x|^d} \,\le\, c_4 r^\alpha |y|^{-d-\alpha}
$$
for some constant $c_4=c_4(d, \alpha)$.
On the other hand, by (\ref{e:int}),
\begin{eqnarray*}
&&\int_{A(0, \sigma, 2r)}
\frac{(\sigma^2-|x|^2)^{\frac{\alpha}2}}{(|y|^2-\sigma^2)^
{\frac{\alpha}2}}\frac1{|y-x|^d}u(y)dy
\,\le\,
c_5 \int_{A(0, \sigma, 2r)}
\frac{r^{\alpha}}{(|y|-\sigma)^
{\frac{\alpha}2}  |y|^{\frac{\alpha}2}      }
\frac1{r^d}u(y)dy\\
&&\le\,
c_5\, r^{\frac{\alpha}2-d}
\int_{A(0, \sigma, 2r)}(|y|-\sigma)^{-\frac{\alpha}2}u(y)dy
\,\le\, c_6r^{-1+\frac{\alpha}2-d}\int_{A(0, \frac{10r}6, 2r)}u(y)
|y|^{-\frac{\alpha}2+1} dy
\end{eqnarray*}
for some constants $c_i=c_i(d, \alpha), i=5, 6$. Hence
\[
\E_x\left[u(Y_{\tau_{D \cap B(z_0, \sigma)}});
Y_{\tau_{D \cap B(z_0, \sigma)}} \in
A(z_0, \sigma, 1- a r) \right]
\,\le\,
c_7\,r^{\alpha}\int_{A(0, \frac{10r}6, 1- a r)}
\frac{u(y)}{|y|^{d+\alpha}}dy
\]
for some constant $c_7=c_7(d, \alpha)$.
\qed

\begin{lemma}\label{l2.2}
Let $D$ be an open set.
Assume that $B(A, \kappa r)\subset D\cap B(Q, r)$
for some $0<r<r_0$ and $\kappa\in (0, \frac12]$.
Suppose that $u\ge0$ is regular harmonic in $D\cap
B(Q, 2r)$ with respect to $Y$
and $u=0$ in $(D^c\cap B(Q, 2r)) \cup B(Q, 1- r)^c$. If $w$ is
a regular harmonic function with respect to $Y$ in $D\cap B(Q, r)$
such that
\[
w(x)=\left\{
\begin{array}{ll}
u(x), & x\in B(Q, \frac{3r}2)^c\cup (D^c\cap B(Q, r))\\
0, & x \in A(Q, r, \frac{3r}2),
\end{array}\right.
\]
then
\[
u(A) \,\ge \,w(A)\, \ge\, C\,\kappa^{\alpha}\,u(x), \quad x \in D
\cap B(Q,\frac32 r)
\]
for some constant $C=C(d, \alpha)>0$.
\end{lemma}

\pf
Without loss of generality, we may assume $Q=0$
and  $x \in D \cap B(0,\frac32 r)$.
The left hand side inequality in the conclusion of the lemma
is obvious, so we only need to prove the right hand side
inequality.
Since $u$ is regular harmonic in $D\cap B(0, 2r)$ with respect
to $Y$ and $u=0$ on
 $B(0, 1- r)^c$, we know from Lemma \ref{l2.U} that
there exists $\sigma\in (\frac{10r}6, \frac{11r}6)$ such that
$$
u(x)=
\E_x\left[u(Y_{\tau_{D \cap B(0, \sigma)}}); \,Y_{\tau_{D \cap
B(0, \sigma)}} \in
A(0, \sigma, 1- r) \right]
\le
c_1r^{\alpha}\int_{A(0, \frac{3r}2, 1- r)}
\frac{u(y)}{|y|^{d+\alpha}}dy
$$
for some constant $c_1=c_1(d, \alpha)$.
On the other hand,
by (\ref{CC1}),  we have  that
\begin{eqnarray*}
&&w(A)\,=\,\int_{A(0, \frac{3r}2 , 1- r )}K^Y_{D\cap B(0, r)}(A, y)u(y)dy
\,\ge\, \int_{A(0, \frac{3r}2, 1- r  )}K^Y_{B(A, \kappa r)}(A, y)u(y)dy\\
&&\ge\, \int_{A(0, \frac{3r}2, 1- r )}K_{B(A, \kappa r)}(A, y)u(y)dy
\,=\,c_2\int_{A(0, \frac{3r}2,  1- r)}\frac{(\kappa r)^{\alpha}}
{(|y-A|^2-(\kappa r)^2)^{\frac{\alpha}2}}\frac1{|y-A|^d}u(y)dy
\end{eqnarray*}
for some constant $c_2=c_2(d, \alpha)$.
Note that $|y-A|\le  2|y|$ on $ A(0, \frac{3r}2, 1- r) $. Hence
\[
w(A)\,\ge\, c_3\,\kappa^{\alpha}\,r^{\alpha}\int_{A(0, \frac{3r}2,  1- r)}
\frac{u(y)}{|y|^
{d+\alpha}}dy
\]
for some constant $c_3=c_3(d, \alpha)$.
Therefore
$
w(A)\,\ge \,c_4 \,\kappa^{\alpha}\,u(x)
$
for some constant $c_4=c_4(d, \alpha)$.
\qed

The following result is a boundary Harnack principle for nonnegative
functions which are harmonic with respect to $Y$ and vanish outside
a small ball. The proof is similar to the
proof in \cite{SW} but we spell out the details for the reader's convenience.

\begin{thm}\label{t2.1}
Suppose that $D$ is  an open set,
$Q\in \partial D$, $r>0$ and that $B(A, \kappa r)$ is a ball in
$D\cap B(Q, r)$.
If $2r < r_0$,
then for any nonnegative functions $u, v$ in $\RR^d$
which are regular harmonic in $D\cap B(Q, 2r)$ with respect to $Y$
and vanish
in $(D^c\cap B(Q, 2r)) \cup B(Q,  1- r)^c$,
we have
$$
C^{-1}\kappa^{d+\alpha}\frac{u(A)}{v(A)}\le
\frac{u(x)}{v(x)}\le C\kappa^{-d-\alpha}\frac{u(A)}{v(A)},
\qquad x\in D\cap B(Q, \frac{r}2),
$$
for some constant $C=C(d, \alpha)>1$.
\end{thm}

\pf
Without loss of generality we may assume that $Q=0$ and $u(A)=v(A)$.
Define $u_1$ and $u_2$ to be regular harmonic functions
in $D\cap B(0, r)$ with respect to $Y$ such that
\[
u_1(x)=\left\{
\begin{array}{ll}
u(x), & r\le |x|<\frac{3r}2\\
0, & x\in B(0, \frac{3r}2)^c\cup(D^c\cap B(0, r))
\end{array}
\right.
\]
and
\[
u_2(x)=\left\{
\begin{array}{ll}
u(x), & x\in B(0, \frac{3r}2)^c\cup(D^c\cap B(0, r))\\
0, & r\le |x|<\frac{3r}2,
\end{array}
\right.
\]
and note that $u=u_1+u_2$. If $D\cap\{r\le |y|<\frac{3r}2\}$
is empty, then $u_1=0$ and the inequality (\ref{e2.6}) below
holds trivially.
So we assume $D\cap\{r\le |y|<\frac{3r}2\}$ is not empty.
Then by Lemma \ref{l2.2},
\[
u(y)\le c_1\kappa^{-\alpha}u(A), \qquad y\in D\cap B(0,
\frac{3r}2),
\]
for some constant $c_1=c_1(d, \alpha)$.
For $x\in D\cap B(0, \frac{r}2)$, we have
\begin{eqnarray*}
u_1(x)&=& \E_x\left[u(Y_{\tau_{D\cap B(0, r)}}): Y_{\tau_{D\cap B(0, r)}}\in
D\cap \{r\le |y|<\frac{3r}2\}\right]\\
&\le&\left(\sup_{D\cap\{r\le |y|<\frac{3r}2\}}u(y)\right)
   \P_x\left( Y_{\tau_{D\cap B(0, r)}}\in
D\cap \{r\le |y|<\frac{3r}2\}\right)              \\
&\le&\left(\sup_{D\cap\{r\le |y|<\frac{3r}2\}}u(y)\right)
   \P_x\left( Y_{\tau_{D\cap B(0, r)}}\in
B(0,r)^c        \right)     \\
&\le&c_1\,\kappa^{-\alpha}\,u(A) \,\P_x\left( Y_{\tau_{D\cap B(0, r)}}\in
B(0,r)^c \right).
\end{eqnarray*}
Now using Lemma \ref{l2.3} we get
\begin{equation}\label{e2.3}
u_1(x)\,\le\, c_2\,\kappa^{-d-\alpha}\,u(A)\,\P_x\left(
Y_{\tau_{(D\cap B(0,r))\setminus
B(A, \frac{\kappa r}2)}} \in B(A, \frac{\kappa r}2)\right),
\qquad x\in D\cap
B(0, \frac{r}2)
\end{equation}
for some constant $c_2=c_2(d, \alpha)$.
Since $2r < r_0$, Corollary \ref{c:Har_1} implies that
\[
u(y)\,\ge\, c_3\,u(A), \qquad y\in B(A, \frac{\kappa r}2)
\]
for some constant $c_3=c_3(d, \alpha)$.
Therefore for $x\in D\cap
B(0, \frac{r}2)$
\begin{equation}\label{e2.4}
u(x) \,=\,  \E_x\left[u(Y_{\tau_{(D\cap B(0, r))\setminus
B(A, \frac{\kappa r}2)}}) \right]     \,\ge\, c_3\,u(A)\,
\P_x\left(Y_{\tau_{(D\cap B(0,r))\setminus
B(A, \frac{\kappa r}2)}} \in B(A, \frac{\kappa r}2)\right).
\end{equation}
From (\ref{e2.3}), the analogue of (\ref{e2.4}) for $v$
and the assumption that $u(A)=v(A)$, it follows that for $x\in D\cap
B(0, \frac{r}2)$,
\begin{equation}\label{e2.6}
u_1(x)\,\le \,c_2\,\kappa^{-d-\alpha}\,v(A)\,
\P_x\left(Y_{\tau_{(D\cap B(0, r)) \setminus
B(A, \frac{\kappa r}2)}} \in B(A, \frac{\kappa r}2)\right)\,\le \,c_4\,
\kappa^{-d-\alpha}\,v(x)
\end{equation}
for some constant $c_4=c_4(d, \alpha).$
Since $u=0$ on $B(0,1- r)^c$, we have that
for $x\in D\cap B(0, r)$,
\begin{eqnarray*}
u_2(x)&=& \int_{A(0, \frac{3r}2,  1-  r)}K^Y_{D\cap B(0, r)}
(x, z)u(z)dz\\
&=&\AA(d, -\alpha)\int_{A(0, \frac{3r}2,  1- r )}
\left(\int_{D\cap B(0, r)}
\frac{G^Y_{D\cap B(0, r)}(x, y)}{|y-z|^{d+\alpha}}dy\right)u(z)dz.
\end{eqnarray*}
Let
\[
s(x)\,:=\,\AA(d, -\alpha)\int_{D\cap B(0, r)}G^Y_{D\cap B(0, r)}(x, y)dy,
\]
then we have
\begin{equation}\label{e2.5}
c_5^{-1}\,\le \,\frac{u_2(x)}{u_2(A)}/\frac{s(x)}{s(A)}\,\le \,c_5,
\end{equation}
for some constant $c_5=c_5(d, \alpha)$.
Applying (\ref{e2.5}) to $u$ and $v$ and Lemma \ref{l2.2} to $v$ and
$v_2$, we obtain for
$x\in D\cap B(0, \frac{r}2)$,
\begin{equation}\label{e2.7}
u_2(x)\,\le\, c_5\,u_2(A)\,\frac{s(x)}{s(A)}\,\le\, c_{5}^2\,
\frac{u_2(A)}{v_2(A)}\,v_2(x)\,
\le\, c_{6}\,\frac{u(A)}{\kappa^{\alpha}v(A)}\,v_2(x)\,=\,
c_{6}\,\kappa^{-\alpha}\,v_2(x),
\end{equation}
for some constant $c_6=c_6(d, \alpha).$
Combining (\ref{e2.6}) and (\ref{e2.7}), we have
\[
u(x)\,\le\, c_{7}\,\kappa^{-d-\alpha}\,v(x), \qquad x\in D\cap B(0,
\frac{r}2),
\]
for some constant $c_{7}=c_{7}(d, \alpha).$
\qed

The theorem above applies to any $\kappa$-fat $D$,
but the harmonic functions there are assumed to vanish outside
a small ball. Thus the theorem above is very useful in studying
properties of positive functions which are harmonic with respect to
$Y$ in $\kappa$-fat sets with diameters less than 1, and not very useful
in the case when the diameters of the $\kappa$-fat sets are large.

Comparing the boundary Harnack principle above with the boundary Harnack
principle for symmetric stable processes in \cite{SW}, we notice that
in the boundary Harnack principle above we assumed an extra
condition that the functions vanish in $B(Q, 1-r)^c$. This extra
condition is not purely technical. In the next section,
we will give an example of
a bounded non-convex domain showing that, without this
extra condition, the boundary Harnack principle for $Y$ fails.

In the  remainder of this section, we will prove a boundary Harnack
principle for nonnegative functions which are harmonic with respect to $Y$
in bounded convex domains without assuming that they vanish
outside small balls.

It is well-known that every convex domain is Lipschitz.
Recall that a bounded domain $D$
is said to be Lipschitz
if there is a localization radius
$R_0>0$  and a constant
$\Lambda >0$
such that
for every $Q\in \partial D$, there is a
Lipschitz  function
$\phi_Q: \R^{d-1}\to \R$ satisfying $\phi_Q (0)= 0$,
$| \phi_Q (x)- \phi_Q (z)| \leq \Lambda
|x-z|$, and an orthonormal coordinate
system $y=(y_1, \cdots, y_{d-1}, y_d):=(\tilde y, y_d)$
such that
$ B(Q, R_0)\cap D=B(Q, R_0)\cap \{ y: y_d > \phi_Q (\tilde y) \}$.
The pair $(R_0, \Lambda)$ is called the
characteristics of the Lipschitz domain $D$. It is easy to see
 that $D$ is $\kappa$-fat with characteristics $(R_0, \kappa_0)$ with some
$\kappa_0=\kappa_0(D)$.

In the remainder of this section we assume $D$ is a bounded
convex domain with the Lipschitz
characteristics  $(R_0, \Lambda)$ and the $\kappa$-fat
characteristics $(R_0, \kappa_0)$.

For every $Q\in \partial D$ and
$ x \in B(Q, R_0)\cap \{ y: y_d > \phi_Q (\tilde y) \}$,
let
$
\delta_Q (x) \,:= \,x_d -  \phi_Q (\tilde x).
$
Since $D$ is bounded Lipschitz, there exists
a constant $c=c(d, \Lambda) \ge 1$ such  for every $Q \in \partial D$
and $ x \in B(Q, R_0)\cap \{ y: y_d > \phi_Q (\tilde y) \}$
 we have
\begin{equation}\label{e:d_com}
c^{-1} \,\delta_Q (x) \,\le\, \rho(x)   \,\le\,    \delta_Q (x)
\end{equation}
The next result is well-known (see Lemma 6.7 in \cite{CZ} for the Brownian
motion case).

\begin{lemma}\label{laa}
 Let $Q \in \partial D$.
Assume that $B(A, \kappa r)\subset D \cap B(Q, r)$
for some positive $r< \frac14 R_0$ and $\kappa\in (0, \frac12]$.
Then there exists $c=c(\alpha, d, D)>0$ such that for every $ y \in
B(Q,(4-\frac12\kappa)r) $ with $\delta_{Q} (y) > \frac12\kappa r$,
we have
$$
  G_{B(Q,4r) \cap D}(A,y) \,\ge\, c \,|A-y|^{-d+\alpha}.
$$
\end{lemma}

\pf
The proof of this result is standard and we omit the details.
\qed

\begin{lemma}\label{l:la}
Let $Q \in \partial D$.
Assume that $B(A, \kappa r)\subset D \cap B(Q, r)$
for some $0<r<\frac14(r_0 \wedge R_0)$ and $\kappa\in (0, \frac12]$.
Suppose that $u\ge0$ is regular harmonic in $D\cap
B(Q, 4r)$ with respect to $Y$
and $u=0$ in $D^c $,
then
\[
u(A)\, \ge\, C\, r^{\alpha}
\int_{A(Q, r, 1+2r)} |z|^{-d-\alpha}u(z)dz
\]
for some constant $C=C(d, \alpha, D)>0$.
\end{lemma}

\pf
Without loss of generality, we may assume $Q=0$.
Let $ \phi (\tilde x):=\phi_0 (\tilde x)$ and $\delta (x):=\delta_0 (x)$.
Since $u$ is regular harmonic in $D\cap
B(0, 4r)$ with respect to $Y$ and $D\cap B(0,r)$ is bounded Lipschitz,
by Theorem \ref{Poisson} we have
\begin{eqnarray*}
&&u(A) = \E_A \left[ u( Y_{\tau_{D \cap B(0,r)}}) \right]
\ge
\int_{A(0,r , 1-r)}K^Y_{D\cap
B(0, r)}(A, z)u(z)dz\\
&&=\int_{A(0, r, 1-r)         }
{\cal A}(d, -\alpha)\int_{D\cap B(0,r ) \cap \{ |y-z| <1\}  }
 \frac{G_{D\cap B(0,r )}^Y(A,y)}{|y-z|^{d+\alpha}}
 dyu(z)dz\\
&&=\int_{A(0, r, 1-r)         }
{\cal A}(d, -\alpha)\int_{D\cap B(0,r ) }
 \frac{G_{D\cap B(0,r )}^Y(A,y)}{|y-z|^{d+\alpha}}
 dyu(z)dz
\end{eqnarray*}
Since $B(A, \kappa r)\subset D \cap B(0, r)$, by the monotonicity of
the Green functions and  (\ref{e:G_0}),
$$
 G_{D\cap B(0,r )}^Y(A,y) \, \ge \,  G_{D\cap B(0,r )}(A,y)  \,
\ge \,G_{B(A, \kappa r)}(A,y),
 \quad y \in B(A, \kappa r).
$$
Thus
$$
u(A) \ge  \int_{A(0,r , 1-r)   }
{\cal A}(d, -\alpha)\int_{ B(A, \kappa r)  }
\frac{G_{B(A, \kappa r)}(A,y)}{|y-z|^{d+\alpha}} dyu(z)dz
= \int_{ A(0, r, 1-r)   }
K_{B(A, \kappa r)}(A, z)u(z)dz,
$$
which is equal to
$$
c_1\int_{A(0, r, 1-r)  }
\frac{(\kappa r)^{\alpha}}{(|z-A|^2-
(\kappa r)^2)^{\frac{\alpha}2}}\frac1{|z-A|^d}u(z)dz
$$
for some constant $c_1=c_1(d, \alpha)$ by (\ref{P_f}).
Note that $|z-A|\le  2|z|$ for $z\in A(0,r , 1-r) $.
Hence
\begin{equation}\label{e:22}
u(A)\,\ge\, c_2\,\kappa^{\alpha}\,r^{\alpha}\int_{A(0,r , 1-r) }
\frac{u(z)}{|z|^
{d+\alpha}}dz
\end{equation}
for some constant $c_2=c_2(d, \alpha)$.
Now we will establish a different lower bound for $u(A)$.
Since $u$ is regular harmonic in $D\cap
B(0, 4r)$ with respect to $Y$ and is zero outside of $D$, we have
\begin{eqnarray*}
u(A) &\ge&
\E_A\left[u(Y_{\tau_{D \cap B(0, 4r)}});
Y_{\tau_{D \cap B(0, 4r)}} \in
A(0, 1-r, 1+2r) \right]\\
&=& \int_{A(0, 1-r, 1+2r) \cap D}K^Y_{B(0, 4r) \cap D}(A, z)u(z)dz
\end{eqnarray*}
Let
$$
\Omega_{\kappa r}\,:=\,\left\{ y \in D  \cap B(0, (4-\frac12\kappa)r)
: \delta (y) > \frac12\kappa r \right\}.
$$
Since
$
  |y-z| \,\le\, |y|+|z|
\,<\,4r +1+2r \,<\,2
$ for $z \in A(0, 1-r, 1+2r)$ and $y \in B(0, 4r)$,
Theorem \ref{Poisson} and (\ref{e:G_0}) imply that
\begin{eqnarray*}
&&K_{B(0,4r) \cap D}^Y(A,z) \,\ge\,
c_3\int_{B(0,4r) \cap D \cap \{ |y-z| <1   \}}
 G^Y_{B(0,4r) \cap D}(A,y)
 dy\\&&\ge\,
c_3\int_{B(0,4r) \cap D \cap \{ |y-z| <1   \}}
 G_{B(0,4r) \cap D}(A,y)
 dy
\, \ge\,
c_3\int_{\Omega_{\kappa r}\cap \{ |y-z| <1  \}}
 G_{B(0,4r) \cap D}(A,y)dy
\end{eqnarray*}
for some constant $c_3=c_3(d, \alpha)$.
By Lemma \ref{laa},  there exists $c_4=c_4(d, \alpha, D)$ such that
$$
  G_{B(0,4r) \cap D}(A,y) \,\ge\, c_4\,
|A-y|^{-d+\alpha}\,\ge\, c_4\,8^{-d+\alpha}
\,r^{-d+\alpha},  \qquad y \in \Omega_{\kappa r}.
$$
So  we have
$$
\inf_{z \in A(0, 1-r, 1+2r) \cap D   }K_{B(0,4r) \cap D}^Y(A,z) \,
\ge\, c_5\,r^{-d+\alpha}
  \inf_{z \in A(0, 1-r, 1+2r) \cap D   }
| \Omega_{\kappa r}     \cap \{ |y-z| <1  \} |
$$
for some constant $c_5=c_5(d, \alpha, D)$.
For each $z \in A(0, 1-r, 1+2r) \cap  D$, let $b^z$
be the point on the line segment
between $z$ and the origin such that
$$
|b^z|\,= \,(3-\frac{\kappa}4)r \quad \mbox{ and } \quad |b^z -z|\,=\,
|z| - (3-\frac{\kappa}4)r.
$$
Note that since $D$ is convex and $z \in D$, $b^z$ is in  $D$.
Let
$$
S_z:=\left\{ (\tilde{y} , y_d) \in B(0,R_0) \cap D:\,\,|\tilde{y}-\tilde{b^z}|
 < \frac{r}{8(1+\Lambda)},\,\,
 \delta(b^z) + \frac12 \kappa r <\delta (y) <  \delta(b^z) + \frac38 r
   \right\}.
$$
We claim that for every $z \in A(0, 1-r, 1+2r) \cap  D$,
$
S_z \subset  \Omega_{ \kappa r}\cap \{ |y-z| <1  \}.
$

For every $y \in S_z$,
\begin{eqnarray*}
&& |y-b^z|
\,\le \,    |\tilde{y}-\tilde{b^z}|+ |y_d-b^z_d|
 \,<\,  \frac{r}{8(1+\Lambda)} +   |  \delta (y)- \delta(b^z)|
+ | \phi (\tilde y)-  \phi (\tilde{b^z})|\\
&&<\,  \frac{r}{8(1+\Lambda)} +  \frac38 r
+ \Lambda |\tilde{y} -  \tilde{b^z}|
\,<\,  \frac{r}{8(1+\Lambda)} +  \frac38 r
+ \frac{r\Lambda}{8(1+\Lambda)}\,= \frac{r}{2}.
\end{eqnarray*}
Thus for every $z \in A(0, 1-r, 1+2r) \cap  D$ and  $y \in S_z$,
 we have
$$
|y-z| \,\le\,  |y-b^z|+ |b^z-z|
\,<\,  \frac{r}{2} + |z|    -
(3-\frac{\kappa}4)r
\,=\,  |z|-2r-\frac14 (2-\kappa)r <1
$$
and
$$
|y| \,\le\, |y-b_z|+|b_z|
\,<\,  \frac{r}{2}+  (3-\frac{\kappa}4)r
\,= \,(\frac{7}{2}-\frac{\kappa}4)r
\,<\, (4-\frac{\kappa}{2})r.
$$
Thus the claim is proved.
Let $\varphi^z( \, \cdot\,):= \phi (\,\cdot\,+\tilde{b^z})-
\phi (\tilde{b^z}) $.
By the change of variable $w=y-b^z$,
 \begin{eqnarray*}
&&|S_z|\,= \,\int_{ \{ |\tilde{y}-\tilde{b^z}| < \frac{r}{8(1+\Lambda)} \} }
 \int_{ \{ b^z_d- \phi (\tilde{b^z}) + \frac12 \kappa r <y_d - \phi
 (\tilde{y}) <b^z_d -  \phi (\tilde{b^z}) + \frac38 r  \}}
dy_d\, d\tilde{y}\\
&&= \int_{\{ |\tilde{w}| < \frac{r}{8(1+\Lambda)}\}}
 \int_{ \{ \frac12 \kappa r <w_d - \varphi^z (\tilde{w}) < \frac38 r \}}
dw_d\,d\tilde{w}
\,\ge\, \int_{ \{ |\tilde{w}| < \frac{r}{8(1+\Lambda)} \}}
 \int_{\{ \frac14 r <w_d - \varphi^z (\tilde{w}) < \frac38 r \} }
dw_d\,d\tilde{w}.
 \end{eqnarray*}
Since  $\varphi^z( \, \cdot\,)$ is Lipschitz  with Lipschitz
constant $\Lambda$
for   $z \in A(0, 1-r, 1+2r) \cap  D$,
the last quantity above is bounded below by $c_6 r^d$ for
some positive constant
$c_6=c_6(D)$ for every $z \in A(0, 1-r, 1+2r) \cap  D$.
Thus we get
$$
\inf_{z \in A(0, 1-r, 1+2r) \cap D   }K_{B(0,4r) \cap D}^Y(A,z) \,\ge\, c_7
\, r^\alpha
$$
for some constant $c_7=c_7(d, \alpha, D)$.
Therefore
\begin{equation}\label{E_3_1}
u(A) \,\ge \,  c_8  \,r^\alpha \int_{
A(0, 1-r, 1+2r) \cap D  }|z|^{-d-\alpha}u(z)dz \,= \,  c_8  \,r^\alpha \int_{
A(0, 1-r, 1+2r)  }|z|^{-d-\alpha}u(z)dz
\end{equation}
for some constant $c_8=c_8(d, \alpha, D)$.
(\ref{e:22}) and (\ref{E_3_1}) imply the lemma.
\qed

The next lemma is a Carleson type estimates for truncated stable processes.

\begin{lemma}\label{l2.2_2_1}
Let $Q \in \partial D$ and assume that $B(A, \kappa r)\subset
D \cap B(Q, r)$
for some $0<r<\frac14(r_0\wedge R_0)$ and $\kappa\in (0, \frac12]$.
If  $u\ge0$ is regular harmonic in $D\cap
B(Q, 4r)$ with respect to $Y$
and $u=0$ in $D^c$,
then
\[
u(A) \,\ge \,C\,u(x), \quad x \in D
\cap B(Q,\frac32 r)
\]
for some constant $C=C(D, d, \alpha, \kappa)$.
\end{lemma}

\pf
Without loss of generality, we may assume $Q=0$.
Let  $x \in D \cap B(0,\frac32 r)$.
Since $u$ is regular harmonic in $D\cap B(0, 4r)$ with
respect to $Y$, we have
\begin{eqnarray*}
u(x)
&=&\E_x\left[u(Y_{\tau_{D \cap B(0, \sigma)}});
Y_{\tau_{D \cap B(0, \sigma)}} \in
A(0, \sigma, 1-\frac{3}{2}r) \right]\\
&&+
\E_x\left[u(Y_{\tau_{D \cap B(0, \sigma)}});
Y_{\tau_{D \cap B(0, \sigma)}} \in
A(0, 1-\frac{3}{2}r, 1+ \sigma) \right]
\,=:\,u_1(x) \,+\, u_2(x),
\end{eqnarray*}
where $\sigma$ is the constant from Lemma \ref{l2.U}.
By Lemma \ref{l2.U},
 We have
\begin{equation}\label{e:11}
u_1(x)\le
c_1r^{\alpha}\int_{   A(0, \frac{10r}6, 1-\frac{3}{2}r) }
\frac{u(z)}{|z|^{d+\alpha}}dz
\end{equation}
for some constant $c_1=c_1(d, \alpha)$.

Now we consider  $u_2(x)$.
We have
$$
u_2(x) =
\E_x\left[u(Y_{\tau_{D \cap B(0, \sigma)}});
Y_{\tau_{D \cap B(0, \sigma)}} \in
A(0, 1-\frac{3}{2}r, 1+\sigma) \right]
\le \int_{A(0, 1-\frac{3}{2}r, 1+\sigma) }K^Y_{B(0, \sigma)}(x, z)u(z)dz.
$$
We know from Lemma \ref{l:P_upper} that
$
K_{B(0, \sigma)}^Y(x,z) \,\le\,c_2 \,r^\alpha
$
for some constant $c_2=c_2(d, \alpha)$.
Therefore
\begin{equation}\label{E_2_1}
u_2(x) \le  c_2 \,r^\alpha \int_{ A(0, 1-\frac{3}{2}r, 1+\sigma)  }u(z)dz\,
\le
\,  c_3 \,r^\alpha \int_{ A(0, 1-\frac{3}{2}r, 1+\sigma)  }
\frac{u(z)}{|z|^{d+\alpha}}  dz,
\end{equation}
for some constant $c_3=c_3(d, \alpha)$ since $\sigma\in (\frac{10r}6,
\frac{11r}6)$.
Combining (\ref{e:11}), (\ref{E_2_1}) and Lemma \ref{l:la},
 we have proved the lemma.
\qed

We shall follow the ``box method'' of \cite{BB}, originally
developed by Bass and Burdzy
 (\cite{BB1} and \cite{BB2}).
Since we are going to use results of \cite{BB},
we will closely follow their notations for the reader's convenience.
Recall that for every $Q\in \partial D$ and
$ x \in B(Q, R_0)\cap \{ y: y_d > \phi_Q (\tilde y) \}$,
$
\delta_Q (x) := x_d -  \phi_Q (\tilde x).
$
For $ x \in B(Q, R_0)\cap \{ y: y_d \ge \phi_Q (\tilde y) \}$, we define
\begin{eqnarray*}
\Delta_Q (x, ar,br) &:=&\{ y: ar >\delta_Q(x) >0,\,
 |\tilde y - \tilde x | < br \} \\
\nabla_Q (x, ar,br) &:=&\{ y: 0 > \delta_Q(x)   >-ar,\,
 |\tilde y - \tilde x | < br \}\\
F^r_{1,Q}&:=& \left\{  X_{\tau_{ \Delta_Q(Q,r,3r)}} \in
\R^d \setminus \left(\Delta_Q(Q,r,3r) \cup \nabla_Q(Q,3r,5r)
\right) \right\}\\
F^r_{2,Q}&:=&
\left\{  X_{\tau_{\Delta_Q(Q,r,3r)}} \in \Delta_Q(Q, 2r,3r) \right\}.
\end{eqnarray*}

We start with the following simple lemma.

\begin{lemma}\label{l:in}
For any positive constants $a$, $b$ and $r$ with $ ((a+b)+b\Lambda)r < R_0$,
$$
 \Delta_Q (Q, ar,br) \,\cup\, \nabla_Q (Q, ar,br) \,\subset\,
B(Q,  ((a+b)+b\Lambda)r)
$$
\end{lemma}

\pf
For $y \in \Delta_Q(Q,ar,br)) \cup \nabla_Q (Q, ar,br)$,
$$
|y-Q| \le |\tilde{y}-\tilde{Q}| + |y_d-Q_d| < br + |\delta_Q (y)|+
|  \phi_Q (\tilde y)
-\phi_Q (\tilde Q)| <  (a+b)r + \Lambda  |\tilde{y}-\tilde{Q}| <
((a+b)+b\Lambda)r.
$$
\qed

Since the above lemma implies that
$$
\Delta_Q (Q, 3r,5r) \,\cup\, \nabla_Q (Q, 3r,5r) \,\subset\,
B(Q,  3(3+2\Lambda)r), $$
the next lemma follows from Lemma 7 and  Remark 2 of \cite{BB} and
the scaling property.

\begin{lemma}\label{l:BB}
Suppose $Q \in \partial D$ and $r <
\frac{R_0}{3(3+ 2 \Lambda)r}$.
If $x \in \Delta_Q (Q, r,3r)$, then $\P_x$ distribution of
$X_{\tau_{\Delta_Q (Q, r,3r)}}$, is
absolutely continuous on $\R^d \setminus ( \Delta_Q (Q, 2r,4r)
\cup \nabla_Q (Q, 2r,4r) )$ with respect to the
$d$-dimensional Lebesgue measure and has a density
function $f^{x,r}_Q$ satisfying
\begin{equation}\label{density}
f^{x, r}_Q(y) \,\le\,  c\, r^{\alpha} \, \P_x(F^r_{1,Q})  \,
(\mbox{\rm dist}(y,\Delta_Q(Q,r,3r))^{-d-\alpha}, \quad y \in
\R^d \setminus ( \Delta_Q (Q, 2r,4r)
\cup \nabla_Q (Q, 2r,4r) )
\end{equation}
where $c=c(D)$ is independent of $Q \in \partial D$.
\end{lemma}

For any Lipschitz function $\psi: \R^{d-1} \to \R$ with Lipschitz constant
$\Lambda$, let
$$
\Delta^\psi  \,:=\,\left\{ y: \frac{R_0}{2(4+3 \Lambda)} >
y_d - \psi (\tilde y) >0,\,
 |\tilde y | < \frac{3R_0}{2(4+3 \Lambda)}  \right\}.
$$
We observe that, for any Lipschitz function $\varphi: \R^{d-1} \to \R$
with Lipschitz constant $\Lambda$,
its dilation $\varphi_r (x) := r \varphi (x/r)$ is also Lipschitz
with the same Lipschitz constant $\Lambda$.
For any $Q \in \partial D$, let $\phi_Q$ be the function in the
definition of a Lipschitz domain.
For any $r>0$, put
$\eta=(2(4+3 \Lambda)r)/R_0$
and $\psi=(\phi_Q)_{\eta}$. Then it is easy to see that
for any $Q \in \partial D$ and $r>0$,
$$
\Delta_Q(Q,r,3r)=\eta \Delta^\psi.
$$
We can show that
$\Delta^\psi \,\subset\,  B(0, \frac12 R_0)$
by the same argument in the proof of Lemma  \ref{l:in}.
On the other hand,
it is easy to see that $\Delta^\psi$
is a $\kappa_1$-fat open set with
$\kappa_1=\kappa_1(\Lambda, R_0)$ for every
Lipschitz function $\psi: \R^{d-1} \to \R$ with Lipschitz constant
$\Lambda$.
Therefore by Proposition \ref{G_2},
there exists positive constant  $r_2 $   such that
for every $Q \in \partial D$ and $r\in (0, r_2]$, we have
\begin{equation}\label{G_4}
 G^Y_{\Delta_Q(Q ,r,3r)}(x, y)
\,\le\, 2\,G_{\Delta_Q(Q ,r,3r)}(x, y),
\quad x, y\in \Delta_Q(Q ,r,3r).
\end{equation}

The next theorem is a boundary Harnack principle for bounded convex
domains and it is the main result of this section. Maybe a word of
caution is in order here. The boundary Harnack principle here is a
little different from the ones proved in \cite{B} and \cite{SW} in
the sense that in the boundary Harnack principle below we require
our harmonic functions to vanish on the whole complement of the
bounded convex domain. However, this will not affect our application
later since we are mainly interested in the case when the harmonic
functions are given by the Green functions of the convex domain.
Recall that $D$ is a bounded convex domain with  Lipschitz
characteristics  $(R_0, \Lambda)$ and $\kappa$-fat characteristics
$(R_0, \kappa_0)$.

\begin{thm}\label{BHP2}
There exist  constants $c >1$  and $r_3 >0 $, depending on
$d, \alpha$ and $D$ such that for any $Q \in \partial D$, $r < r_3$ and any
nonnegative  functions $u, v$  which are regular harmonic
with respect to $Y$ in $D \cap B(Q, 6(3+ 2\Lambda ) r)$ and vanish in
$D^c$, we have
\begin{equation}\label{e:BHP}
\frac{ u(x)}{ v(x)}
\, \le \,c\, \frac{ u(y)}{ v(y)}  \quad  \mbox{ for any }
x,y \in   D  \cap B(Q, \frac{r}{1+ \Lambda}).
\end{equation}
\end{thm}

\pf
Without loss of generality, we may assume $Q=0$.
Let $r_3 :=  (r_0 \wedge r_2 \wedge R_0)/(6(3+ 2\Lambda ))$
and fix a $r < r_3$.
For notational convenience, we denote
 $ \phi (\tilde x):=\phi_0 (\tilde x)$,  $\delta (x):=\delta_0 (x)$,
$F_1=F^r_{1,0}$, $F_2:=F^r_{2,0}$,
$\Delta_0(0 ,a,b):= \Delta(a,b)$ and $\Delta:=\Delta(r,3r)$.
Note that if $y \in  B(0, \frac{r}{1+ \Lambda}) \cap D $, then
$$
y_d -\phi (\tilde y) \le  y_d +  \Lambda |\tilde y| < \frac{r}{1+ \Lambda}
 +  \Lambda \frac{r}{1+ \Lambda} =r.
$$
So $B(0, \frac{r}{1+ \Lambda}) \cap D \subset  \Delta(r,r)$.
Thus it is enough to
 consider $x, y \in  \Delta(r,r) $.
By Lemma \ref{l:in},
 $ \Delta(2r,4r) \subset  B(0,2(3+ 2\Lambda )r  )$.
Thus by Lemma \ref{l2.2_2_1}, there exists a positive constant
$c_1=c_1(D, d, \alpha)$
\begin{equation}\label{partial^s}
\sup_{\Delta(2r,4r)} u \,\le\, c_1\, u(A)
\end{equation}
where
 $B(A, \kappa_0 r ) \subset D \cap  B(0, \frac43(3+ 2\Lambda )r) $.
Let $\Delta_1:= \{ z \in D \setminus \Delta: \mbox{dist}(z, \Delta)<1\}$.
By  (\ref{G_4}) and Theorem \ref{Poisson},
we have
\begin{equation}\label{e:Del}
K_{\Delta}^Y(x,z) \,\le\, 2\, K_{\Delta}(x,z), \quad z \in  \Delta_1.
\end{equation}
Thus,
\begin{eqnarray*}
\E_x\left[ u(Y_{\tau_{\Delta}}):
 Y_{\tau_{\Delta}} \in  D \setminus\Delta( 2r,4r) \right]
&=&
\int_{  \Delta_1 \setminus\Delta( 2r,4r)}
K_{\Delta}^Y(x,z) u(z)dz\\
&\le&
2\,\E_x\left[ u(X_{\tau_{\Delta}}):
X_{\tau_{\Delta}} \in \Delta_1 \setminus\Delta( 2r,4r) \right]
.
\end{eqnarray*}
By Lemma \ref{l:BB}, we have
\begin{eqnarray*}
\E_x\left[ u(X_{\tau_{\Delta}}):
 X_{\tau_{\Delta}} \in \Delta_1 \setminus\Delta(0, 2r,4r) \right]
&\le& c_2\,  r^{\alpha} \,\P_x(F_1)   \,
\int_{\Delta_1 \setminus\Delta(0, 2r,4r)}
     (\mbox{dist}(z,\Delta))^{-d-\alpha}       u(z)dz\\
&\le& \,c_3\,  r^{\alpha} \,\P_x(F_1)  \,
\int_{\Delta_1 \setminus\Delta(0, 2r,4r)}
     |z|^{-d-\alpha}     u(z)dz
\end{eqnarray*}
for some constants $c_i=c_i(d, \alpha, D)$, $i=2, 3$.
In the last inequality above, we have used the fact that for $z \in
\Delta_1 \setminus\Delta( 2r,4r)$,
$$
|z|\, \le \,\mbox{dist}(z,\Delta) + \mbox{diam}(\Delta) \, \le \,
 \mbox{dist}(z,\Delta) +2(4+ 3\Lambda )r
 \,\le\,  3(3+ 2\Lambda )\, \mbox{dist}(z,\Delta).
$$
Therefore, by (\ref{partial^s}) and (\ref{e:Del}),
for every $x \in \Delta(r,r)$,
\begin{eqnarray}
u(x) &=& \E_x\left[ u(Y_{\tau_{\Delta}}):  Y_{\tau_{\Delta}}
\in   \Delta( 2r,4r)      \right]
 \,+\,  \E_x\left[ u(Y_{\tau_{\Delta}}):  Y_{\tau_{\Delta}}
\in  D \setminus\Delta( 2r,4r) \right]\nonumber\\
 &\le& c_1\, u (A)\, \P_x\left(  Y_{\tau_{\Delta}} \in
\Delta( 2r,4r) \right)
 \,+\, 2\,\E_x\left[ u(X_{\tau_{\Delta}}):  X_{\tau_{\Delta}}
\in  \Delta_1 \setminus\Delta( 2r,4r) \right]\nonumber\\
 &\le& 2\, c_1\, u (A)\, \P_x\left(  X_{\tau_{\Delta}} \in
\Delta( 2r,4r) \right)
 \,+\, 2\,c_3\, r^{\alpha} \,\P_x(F_1)      \int_{\Delta_1
\setminus\Delta( 2r,4r)}
     |z|^{-d-\alpha}       u(z)dz\nonumber \\
&=& 2\,c_1 \, u (A)\,\P_x(F_1)
 \,+\, 2\,c_3\, r^{\alpha} \,\P_x(F_1)      \int_{\Delta_1
\setminus\Delta( 2r,4r)}
     |z|^{-d-\alpha}       u(z)dz. \label{e:3333}
\end{eqnarray}
On the other hand, by Lemma \ref{l:la},
\begin{equation}\label{BHP11}
u(A) \,\ge\, c_4 \, r^{\alpha}
\int_{   A(0, \frac{9}{4}(3+ 2\Lambda ) r, 1+3(3+ 2\Lambda ) r) }
     |z|^{-d-\alpha}       u(z)dz
\end{equation}
for some constant $c_4=c_4(D, d, \alpha)$.
Moreover, since $u$ is harmonic with respect to $Y$ in
$B(A, \frac12 \kappa_0 r)$
and $B(0,\frac{9}{4}(3+ 2\Lambda ) r) \subset B(A, 1-\frac12 \kappa_0 r)$,
 by (\ref{CC1})
\begin{eqnarray*}
&&u(A) \,\ge
\, \E_A\left[ u(Y_{\tau_{B(A, \frac12\kappa_0  r)}}):
Y_{\tau_{B(A, \frac12 \kappa_0   r)}}
\in  B(0,\frac{9}{4}(3+ 2\Lambda ) r) \setminus    \Delta( 2r,4r)   \right]\\
& &= \int_{B(0,\frac{9}{4}(3+ 2\Lambda ) r) \setminus  \Delta( 2r,4r) }
K^Y_{B(A, \frac12\kappa_0  r)} (A,z) u(z)dz
 \ge \int_{B(0,\frac{9}{4}(3+ 2\Lambda ) r) \setminus  \Delta( 2r,4r) }
K_{B(A, \frac12\kappa_0  r)} (A,z) u(z)dz.
\end{eqnarray*}
So by (\ref{P_f}),
\begin{eqnarray}
u(A)   &\ge&
c_5\int_{B(0,\frac{9}{4}(3+ 2\Lambda ) r) \setminus  \Delta( 2r,4r) }
\frac{(\frac12 \kappa_0 r)^{\alpha}}{(|z-A|^2-
(\frac12 \kappa_0 r)^2)^{\frac{\alpha}2}}\frac1{|z-A|^d}u(z)dz\nonumber\\
&\ge& c_6 r^{\alpha}
\int_{B(0, \frac{9}{4}(3+ 2\Lambda ) r) \setminus  \Delta( 2r,4r)  }
|z|^{-d-\alpha}   u(z)dz \label{BHP12}
\end{eqnarray}
for some constant $c_i=c_i(D, d, \alpha)$, $i=5, 6$.
Since
$\Delta_1 \subset B(0, 1+3(3+ 2\Lambda )r)$,
by combining (\ref{BHP11}) and (\ref{BHP12})
 we get
\begin{equation}\label{e:333}
u(A) \ge c_7 r^{\alpha}
\int_{B(0, 1+3(3+ 2\Lambda )r) \setminus \Delta( 2r,4r)  } |z|^{-d-\alpha}
u(z)dz
 \ge c_7 r^{\alpha}
 \int_{\Delta_1 \setminus\Delta( 2r,4r)}
     |y|^{-d-\alpha}       u(z)dz
\end{equation}
for some constant $c_7=c_7(D, d, \alpha)$.
Putting (\ref{e:3333}) and (\ref{e:333}) together, we have
\begin{equation}\label{e:BHP0}
 u(x)\, \le \,c_8\, u(A)\, \P_x(F_1) \,=\, c_8\frac{u (A)}{v(A)} \,
v(A)\,\P_x(F_1)
\end{equation}
for some constant $c_8=c_8(D, d, \alpha)$.
By Lemma 6 in \cite{BB}, we have
\begin{eqnarray}
&&\P_x(F_1) \,\le\,
c_9\,\P_x(F_2)
=  c_9\, \P_x( X_{\tau_{\Delta}} \in \Delta( 2r,3r))
=c_9 \int_{\Delta( 2r,3r) }
K_{\Delta}(x,z)dz\nonumber\\
&&=c_9  \int_{\Delta( 2r,3r) }  {\cal A}(d, -\alpha)\int_{  \Delta   }
\frac{G_\Delta(x,y)}{|y-z|^{d+\alpha}} dydz
\le  c_9\int_{\Delta( 2r,3r) }  {\cal A}(d, -\alpha)\int_{  \Delta   }
\frac{G_\Delta^Y(x,y)}{|y-z|^{d+\alpha}}
dydz \label{e:yes}
\end{eqnarray}
for some constant $c_9=c_9(D, d, \alpha)$.
We have used (\ref{e:G_0}) in the last inequality above.
For every $z \in \Delta( 2r,3r)$ and $y \in  \Delta $,
$
|y-z| \le |y|+|z| <4(3+ 2\Lambda )r <r_0 <1.
$
So (\ref{e:yes}) is, in fact, equal to
$$
c_9  \int_{\Delta( 2r,3r) }
K^Y_{\Delta}(x,z)dz  \,=\,c_9\, \P_x( Y_{\tau_{\Delta}} \in \Delta( 2r,3r) ).
$$
By the Harnack inequality (Theorem \ref{T:Har}),
\begin{equation}\label{e:BHP1}
 v(A)\,\P_x( Y_{\tau_{\Delta}} \in \Delta( 2r,3r) )
\,\le\, c_{10}\,\E_x[ v(Y_{\tau_{\Delta}}): Y_{\tau_{\Delta}}
\in \Delta( 2r,3r) ]
\,\le \,c_{10}\, v(x)
\end{equation}
for some constant $c_{10}=c_{10}(D, d, \alpha)$.
From (\ref{e:BHP0}) and (\ref{e:BHP1}), we conclude
\begin{equation}\label{e:BHP2}
\frac{ u(x)}{ v(x)}
\, \le \,c_{11}\, \frac{ u(A)}{ v(A)}, \quad x \in  \Delta(r,r)
\end{equation}
for some constant $c_{11}=c_{11}(D, d, \alpha)$.
The above argument also implies that
$$
\frac{ v(y)}{ u(y)}
\, \le \,c_{11}\, \frac{ v(A)}{ u(A)},  \quad y \in  \Delta(r,r).
$$
Therefore
$$
\frac{ u(x)}{ v(x)} \, \le
\,c_{11}\, \frac{ u(A)}{ v(A)}  \, \le
\,c_{11}^2\,
\frac{ u(y)}{ v(y)}, \quad x, y \in  \Delta(r,r).
$$
\qed

In the remainder of this section, we fix $r_3 >0$ from Theorem \ref{BHP2}.
The following result is analogous to Lemma 5 of \cite{B}. We recall from
Definition \ref{fat} that  for each
$z \in \partial D$ and $r \in (0, R_0)$, $A_r(z)$ is a point
in $D \cap B(z,r)$
satisfying $B(A_r(z),\kappa_0 r)  \subset D \cap B(z,r)$.

\begin{lemma}\label{l:5B}
There exist positive constants $C=C(D, d, \alpha)$ and
$\gamma=\gamma(d, \alpha)<
\alpha$ such that for any $Q\in \partial D$ and
$r\in (0, r_3)$, and nonnegative function
$u$ which is harmonic with respect to $Y$ in $D \cap B(Q, r)$ we have
\[
u(A_s(Q))\ge C(s/r)^{\gamma}u(A_r(Q)), \qquad s\in (0, r).
\]
\end{lemma}

\pf
Without loss of generality, we may assume $Q=0$.
Fix $r < r_3$ and
let
$$
r_k\,:=\,\left(\frac{2}{\kappa_0}\right)^{-k}r,
 \quad A_k\,:=\, A_{r_k}(0)
\quad \mbox{ and } \quad B_k\,:=\,B(A_k, r_{k+1}), \quad k=0,1, \cdots.
$$
Note that the $B_k$'s are disjoint. So by the harmonicity of $u$, we have
$$
u(A_k)
\,\ge\, \sum_{l=0}^{k-1} \E_{A_k}\left[u(Y_{\tau_{B_k}}):\,
Y_{\tau_{B_k}} \in B_l \right]\\
\,=\, \sum_{l=0}^{k-1} \int_{B_l} K_{B_k}^Y(A_k, z) u(z) dz.
$$
Since $r < r_3$,  (\ref{CC1}) and Corollary \ref{c:Har_1} imply that
$$
 \int_{B_l} K_{B_k}^Y(A_k, z) u(z) dz \,\ge\, c_1\, u(A_l)
\int_{B_l} K_{B_k}(A_k, z) dz
$$
for some constant $c_1=c_1(d, \alpha)$.
Using  the explicit formula of  $K_{B_k}$, one can easily check that
$$
\int_{B_l}K_{B_k}(A_k, z)dz
\,\ge\, c_2 \, \left(\frac{\kappa_0}{2}\right)^{-(k-l) \alpha},
\quad z \in B_l,
$$
for some constant $c_2=c_2(d, \alpha)$.
Therefore,
$$
\left(\frac{2}{\kappa_0}\right)^{k\alpha}
 u(A_k) \,\ge \, c_3 \sum_{l=0}^{k-1} \left(\frac{2}{\kappa_0}
\right)^{l\alpha} u(A_l)
$$
for some constant $c_3=c_3(d, \alpha)$.
The remainder of the proof is  same as in the proof of
Lemma 5 in \cite{B} and so we omit it.
\qed

The next lemma is analogous to Lemma 14 of \cite{B}.

\begin{lemma}\label{l:14B}
Suppose $M:= 6(3+ 2\Lambda )(1+\Lambda ) $.
Then there exist positive constants $c_1=c_1(D,d,\alpha)$  and
$c_2=c_2(D,d,\alpha)<1$ such that
for any $Q \in \partial D $, $r < \frac{r_3}{1+\Lambda}$ and nonnegative
function $u$ which is regular harmonic with respect to $Y$
in $D \cap B(Q, M r)$ and vanishes in  $D^c$,
$$
\E_x\left[u(Y_{\tau_{D \cap B_k}}):\,Y_{\tau_{D
\cap B_k}} \in A(Q, r, 1+M^{-k}r)
 \right] \,\le\, c_1\,c_2^{k} \, u(x), \quad x \in D \cap B_k,
$$
where $B_k:=B(Q, M^{-k}r)$, $ k=0,1, \cdots$.
\end{lemma}

\pf
Without loss of generality, we may assume $Q=0$. Fix
$r <\frac{r_3}{1+\Lambda}$ and
a nonnegative function $u$ which is harmonic with respect to $Y$
in $D \cap B(0, Mr)$ and  vanishes in  $\R^{d}\setminus D$.

Let $r_k:=M^{-k}r $, $B_k:=B(0,r_k)$ and
$$
u_k(x)  \,:=\, \E_x\left[u(Y_{\tau_{D \cap B_k}}):
Y_{\tau_{D \cap B_k}}
 \in A(0, r, 1+r_k) \right], \quad x \in D \cap B_k.
$$
Note that
\begin{eqnarray*}
 u_{k+1}(x) &=&  \E_x\left[u(Y_{\tau_{D \cap B_{k+1}}}):\,
Y_{\tau_{D \cap B_{k+1}}}
 \in A(0, r, 1+r_{k+1}) \right] \\
 &=&  \E_x\left[u(Y_{\tau_{D \cap B_{k+1}}}):\,
\tau_{D \cap B_{k+1}} = \tau_{D \cap B_{k}}, ~
Y_{\tau_{D \cap B_{k+1}}}
 \in A(0, r, 1+r_{k+1})\right] \\
&=&  \E_x\left[u(Y_{\tau_{D \cap B_{k}}}):\,
\tau_{D \cap B_{k+1}} = \tau_{D\cap B_{k}}, ~
Y_{\tau_{D \cap B_{k}}}
 \in A(0, r, 1+r_{k+1}) \right] \\
&\le & \E_x\left[u(Y_{\tau_{D \cap B_{k}}}):\,
Y_{\tau_{D \cap B_{k}}}
 \in A(0, r, 1+r_{k+1}) \right]\\
&\le & \E_x\left[u(Y_{\tau_{D \cap B_{k}}}):\,
Y_{\tau_{D \cap B_{k}}}
 \in A(0, r, 1+r_{k}) \right]
\end{eqnarray*}
Thus
\begin{equation}\label{e:dec1}
u_{k+1}(x) \,\le\, u_{k}(x).
\end{equation}
Let
$A_k\,:=\,A_{r_k}(0) $.
We have
 \begin{eqnarray*}
 u_{k}(A_k) &=&  \E_{A_k}\left[u(Y_{\tau_{D\cap B_k}}):\,
Y_{\tau_{D \cap B_{k}}}
 \in A(0, r, 1+r_{k}) \right] \\
&\le&  \E_{A_k}\left[u(Y_{\tau_{ B_{k}}}):\,
Y_{\tau_{ B_{k}}}
 \in A(0, r, 1+r_{k})\right]\\
&\le&  \int_{A(0, r, 1-r_{k})} K_{B_k}^Y(A_k,z) u(z) dz
\,+\,  \int_{A(0, 1-r_k, 1+r_{k})} K_{B_k}^Y(A_k,z) u(z) dz.
\end{eqnarray*}
For $z \in A(0, r, 1-r_{k})$, by (\ref{CC2}) and (\ref{P_f})
 we get
$$
K_{B_k}^Y(A_k,z) \,\le\,2\,K_{B_k} (A_k,z) \,\le\,  c_1
\,\frac{M^{-k\alpha} r^\alpha}{|z|^{d+\alpha}}
$$
for some constant $c_1=c_1(d, \alpha)$.
For $z \in A(0,1-r_k, 1+r_{k})$, we use Lemma \ref{l:P_upper} and  get
 $$
K_{B_k}^Y(A_k,z) \,\le\, c_2\,
\,M^{-k\alpha} r^\alpha
\,\le\, c_3\,
\,\frac{M^{-k\alpha} r^\alpha}{|z|^{d+\alpha}}
$$
for some constant $c_i=c_i(d, \alpha)$, $i=2, 3$.
Therefore
\begin{equation}\label{e:dec2}
 u_{k}(A_k) \,\le\, c_4\,M^{-k\alpha}  r^\alpha\int_{A(0, r, 1+r_{k})}
u(z)\frac{dz}{|z|^{d+\alpha}}
\end{equation}
for some constant $c_4=c_4(d, \alpha)$.
From Lemma \ref{l:la}, we have
\begin{equation}\label{e:dec3}
 u(A_0) \,\ge\, c_5\ r^\alpha\int_{A(0, r, 1+r)}
u(z)\frac{dz}{|z|^{d+\alpha}}
\end{equation}
for some constant $c_5=c_5(D, d, \alpha)$.
(\ref{e:dec2}) and (\ref{e:dec3}) imply that
$ u_{k}(A_k) \,\le\, c_6\,M^{-k\alpha}  u(A_0)$
for some constant $c_6=c_6(D, d, \alpha)$.
On the other hand, using Lemma \ref{l:5B}, we get
$u(A_0) \,\le\, c_7\,M^{k\gamma } u(A_k)$
for some constant $c_7=c_7(D, d, \alpha)$.
Thus, $ u_{k}(A_k) \,\le\, c_6c_7\,M^{-k(\alpha-\gamma)}u(A_k)$.
By (\ref{e:dec1}) and (\ref{e:BHP2}), we have
$$
\frac{u_k(x)}{u(x)}    \,\le\,
\frac{u_{k-1}(x)}{u(x)}
\,\le\, c_8\, \frac{u_{k-1}(A_{k-1})}{u(A_{k-1})}  \,\le\,
c_6c_7c_8\,M^{-k(\alpha-\gamma)}
$$
for some constant $c_8=c_8(D, d, \alpha)$.
\qed

Now the next theorem follows from Lemma \ref{l:5B}, Theorem \ref{BHP2}
and Lemma \ref{l:14B}
(instead of using Lemma 5,  Lemma 13 and Lemma 14 in \cite{B} respectively)
in very much the same way as in the case of symmetric stable process proved
in Lemma 16 of \cite{B}. We omit the details.

\begin{thm}\label{t2.2}
There exist positive constants $r_4$, $M_1$,
$C$ and $\nu$ depending on $D$ and $\alpha$ such
that for any $Q \in \partial D $, $r < r_4$ and nonnegative
functions $u, v$ which are regular harmonic with respect to $Y$
in $D \cap B(Q, M_1 r)$, vanish in  $\R^{d} \setminus D$, and
satisfy $u(A_r(Q))=v(A_r(Q))>0$, we have
\[
\left|\frac{u(x)}{v(x)}-\frac{u(y)}{v(y)}\right|
\le C\left(\frac{|x-y|}r\right)^{\nu},
\qquad x, y\in  D \cap B(Q, r).
\]
In particular, the limit $\lim_{D \ni x\to w}u(x)/v(x)$ exists
for every $w\in \partial  D \cap B(Q, r)$.
\end{thm}

Using the results above and repeating the arguments in the proof of
Theorem 4.1 of \cite{SW} we can get the following result
identifying the Martin boundary of any
bounded convex domain. For the definition and basic results on
Martin boundary, one can see \cite{KW} and \cite{SW}.

\begin{thm}\label{T:mart2}
Suppose that $D$ is a bounded convex domain in $\R^d$.
Then both the Martin boundary and the minimal Martin boundary of $D$
with respect to $Y$ coincide with the Euclidean boundary
of $D$.
\end{thm}

\pf We omit the details.
\qed

\section{Counterexample}

In this section, we present an example of a bounded non-convex domain
on which the boundary Harnack principle for $Y$ fails.

Consider the domain in $R^d$
$$
D:=(-100, 100)^{d} \setminus \left((-100, 50]^{d-1} \times
 [-1/2, 0]\right).
$$
Of course there are nothing special about the numbers 100 and 50
above, they are just two big numbers.

Suppose the boundary Harnack principle (not necessarily scale
invariant) is true for $D$ at the origin. i.e., there exist
constants $R_1 >0 $ and  $M_1 >1$ such that for any  $r < R_1$ and
any nonnegative  functions $u, v$  which are regular harmonic with
respect to $Y$ in $D \cap B(0, M_1 r)$ and vanish in $D^c$, we have
\begin{equation}\label{ce0}
\frac{ u(x)}{ v(x)}
\, \le \,c\, \frac{ u(y)}{ v(y)}  \quad  \mbox{ for any }
 x,y \in
D  \cap B(0, r),
\end{equation}
where $c=c(D,r)>0$ is independent of harmonic functions $u$ and $v$.
Choose an $r_1 < R_1$ with $M_1 r_1 <1/2$ and let $A:=(\tilde 0,
\frac12 r_1)$. We define a function $v$ by
$$
v(x):=\P_x \left(Y_{ \tau_{D \cap B(0,M_1r_1)}} \in \{y \in D; y_d
>0\}\right).
$$
By definition $v$ is regular harmonic in $D \cap B(0,M_1r_1)$ with
respect to $Y$ and vanishes in $D^c$. Applying $v$ above to
(\ref{ce0}), we have a Carleson type estimate at $0$, i.e., there
exists  constant $c_1=c_1(D, r_1)>0$ such that for any nonnegative
function $u$  which is regular harmonic with respect to $Y$ in $D
\cap B(0, M_1 r_1)$ and vanishes in $D^c$ we have
\begin{equation}\label{ce1}
u(A) \,\ge \,c_1\,u(x), \quad x \in D \cap B(0, r_1).
\end{equation}
We will construct a bounded positive function $u$ which is regular
in $D \cap B(0, M_1 r)$ with respect to $Y$ and vanishes
in $D^c$ for which (\ref{ce1}) fails.

For $n \ge 1$, we put
\begin{eqnarray*}
C_n&:=&\left\{ (\tilde x,x_d) \in D;\quad |\tilde x| \le \frac{r_1}8,
\quad  x_d
\le -1+2^{-n}r_1^2 \right\}\\
D_n&:=&\left\{ (\tilde y,y_d) \in D;\quad y_d>0,  \quad
|x-y| <1 \quad \mbox{ for some } x \in  C_n\right\}.
\end{eqnarray*}
Note that $D_n \supset D_{n+1} \supset \cdots$ and
$\cap_{1}^{\infty} D_n= \emptyset$.
Moreover, it is easy to see that
\begin{equation}\label{ce2}
D_n \subset B(0,r_1) \cap D,  \quad \mbox{ for } n \ge 3.
\end{equation}
In fact, for any $y\in D_n$, we have $y_d\in (0, 2^{-n}r_1^2)$ and
$|y-x |<1$ for some $x \in C_n$, thus $y_d-x_d > -x_d
\ge 1-2^{-n}r_1^2$ and
$$
 |\tilde y - \tilde x|^2 + |y_d  - x_d|^2 \,=\,|(\tilde y,y_d)-
(\tilde x, y_d)|^2
+ |(\tilde x,y_d)-(\tilde x, x_d)|^2 \,<\,1.
$$
Hence
$$
|\tilde y| \,\le\, |\tilde x|+ |\tilde y - \tilde x|
\,\le\, \frac{r_1}8 + \sqrt{ 1- |y_d  - x_d|^2} \,<\, \frac{r_1}8+
\sqrt{ 2^{-n+1}r_1^2 }.
$$
Since $r_1 <1$, we get for $n \ge 3$
$$
|\tilde y|^2 + y_d^2 \,< \,2^{-n}r_1^2 + (\frac{r_1}8+
2^{(-n+1)/2}r_1 )^2 \,<\, r_1^2.$$

For any $n$, let $T_{D_n}$ be the first hitting time of $D_n$ by the
process $Y$. Note that since $\cap_{1}^{\infty} D_n= \emptyset$,
$$
\P_A(\tau_{D \cap B(0,M_1r_1)} > T_{D_n}) \,\to\, 0, \quad
\mbox{ as } n \to \infty.
$$
Choose $n\ge 3$ large so that
\begin{equation}\label{ce3}
\P_A(\tau_{D \cap B(0,M_1r_1)} > T_{D_{n}}) \, < \, \frac{c_1}{2}
 \end{equation}
and define
$$
u(x)\,:=\P_x \left( Y_{ \tau_{D \cap B(0,M_1r_1)}} \in C_{n}\right).
$$
$u$ is a nonnegative bounded function which is
regular harmonic in $D \cap B(0,M_1r_1)$ with respect to
$Y$ and vanishes in $D^c$. It also vanishes
continuously on $\partial D \cap B(0,M_1r_1)$.
Note that by Theorem \ref{Poisson},
$$\P_{A} \left( Y_{ \tau_{D \cap B(0,M_1r_1)}} \in C_{n},
\, \tau_{D \cap B(0,M_1r_1)} \le T_{D_{n}}   \right)
\, =\,
\P_{A} \left(    Y_{ \tau_{D \cap B(0,M_1r_1)
\setminus D_{n}}} \in C_{n}      \right)=0.
$$
Thus by the strong Markov property,
\begin{eqnarray*}
u(A) &=&\P_A \left( Y_{ \tau_{D \cap B(0,M_1r_1)}} \in C_{n},\,
\tau_{D \cap B(0,M_1r_1)} > T_{D_{n}} \right)\\
&=&\E_A \left[\P_{Y_{T_{D_{n}}}}
\left( Y_{ \tau_{D \cap B(0,M_1r_1)}} \in C_{n} \right);\,
\tau_{D \cap B(0,M_1r_1)} > T_{D_{n}}   \right]\\
&\le&\P_A\left(\tau_{D \cap B(0,M_1r_1)} > T_{D_{n}}\right) \left(
\sup_{x \in D_{n}} u(x)\right) \, < \, \frac{c_1}2 \left(\sup_{x \in
D \cap B(0,r_1)   } u(x)\right).
\end{eqnarray*}
In the last inequality above, we have used (\ref{ce2})-(\ref{ce3}).
But by (\ref{ce1}), $u(A) \,\ge \,c_1   \,  \sup_{x \in D \cap
B(0,r_1)   }  u(x)$, which gives a contradiction. Thus the boundary
Harnack principle is not true for $D$ at the origin.

By smoothing off the corners of $D$, we can easily construct
a smooth bounded non-convex domain on which the boundary Harnack principle
fails for the truncated stable process $Y$.

\vspace{.3cm} \noindent {\bf Acknowledgment}: We thank Zoran
Vondracek for helpful comments. We also thank an anonymous referee
for helpful comments on the first version of this paper.

\vspace{.1in}
\begin{singlespace}

\end{singlespace}
\end{doublespace}

\end{document}